\documentclass[a4paper,10pt]{amsart}
\usepackage[foot]{amsaddr}
\usepackage{geometry}
\usepackage{amsfonts}
\usepackage{amsthm}
\usepackage{amssymb}
\usepackage{amsmath}
\usepackage{mathrsfs}
\usepackage{comment}
\usepackage{theoremref}
\usepackage{enumerate}
\usepackage{bbm}
\usepackage{bm}
\usepackage{comment}
\usepackage{mathtools}
\usepackage{a4wide}
\usepackage{xcolor}
\usepackage{lipsum}

\makeatletter
\newcommand{\proofpart}[2]{%
    \par
  \addvspace{\medskipamount}%
  \noindent\emph{Step #1: #2}\par\nobreak
  \addvspace{\smallskipamount}%
  \@afterheading
}
\makeatother

\DeclarePairedDelimiter\abs{\lvert}{\rvert}%
\DeclarePairedDelimiter\norm{\lVert}{\rVert}%

\makeatletter
\let\oldabs\abs
\def\abs{\@ifstar{\oldabs}{\oldabs*}}
\let\oldnorm\norm
\def\norm{\@ifstar{\oldnorm}{\oldnorm*}}
\makeatother

\makeatletter
\g@addto@macro\bfseries{\boldmath}
\makeatother

\newcommand{\T}{\mathbb{T}}

\newcommand{\Ka}{\mathcal{K}}
\newcommand{\conj}[1]{\overline{#1}}
\newcommand{\D}{\mathbb{D}}

\newcommand{\dist}[2]{\text{dist}( #1, #2 ) }

\renewcommand\Re{\operatorname{Re}}

\newcommand{\supp}[1]{\text{supp}({#1})}

\newtheorem{thm}{Theorem}[section]
\newtheorem{lemma}[thm]{Lemma}

\newtheorem{cor}[thm]{Corollary}
\newtheorem{prop}[thm]{Proposition}

\theoremstyle{definition}

\theoremstyle{definition}

\newcommand{\Addresses}{{
		\bigskip
		\footnotesize
		
		Adem Limani, \\ \textsc{Centre for Mathematical Sciences \\ Lund University \\
		Lund, Sweden}\\
		\texttt{adem.limani@math.lu.se}
		
	}}

\begin{document}
\title{\textbf{Generic measures with slowly decaying Fourier coefficients}} 

\author{Adem Limani} 
\address{Centre for Mathematical Sciences, Lund University, Sweden}
\email{adem.limani@math.lu.se}

\date{\today}

\begin{abstract}
We investigate threshold phenomena in weighted $\ell^2$-spaces and characterize the critical regimes where elements with either small support or maximally bad range can be constructed. Our results are shown to be optimal in several respects, and our proofs principally rely on techniques involving sparse Fourier spectrum. We further show that these seemingly pathological constructions are actually generic from certain categorical perspectives. 
\end{abstract}

\maketitle
\section{Introduction}\label{SEC:INTRO}
\subsection{Supports of elements in weighted $\ell^2$}
Let $w=(w_n)_{n=0}^\infty$ be positive numbers and denote by $\ell^2(w)$ the Hilbert space of distributions $S$ on $\T$ satisfying 
\[
\norm{S}_{\ell^2(w)} = \left( \sum_{n} \abs{\widehat{S}(n)}^2 w_{|n|} \right)^{1/2}< \infty,
\]
where $\widehat{S}(n) := S(\zeta^{-n} )$ denotes the $n$-th Fourier coefficient of $S$. Since any distribution $S$ on $\T$ is uniquely determined by its Fourier coefficients $(\widehat{S}(n))_n$, we allow a slight abuse of notation and write $S\in \ell^2(w)$ interchangeably with $(\widehat{S}(n))_n \in \ell^2(w)$. The classical $\ell^2$-space is obtained when $(w_n)_n$ is a constant sequence, which by Parseval's Theorem can be identified with $L^2(\T,dm)$, the space of square integrable functions $f$ on the unit-circle $\T$:
\[
\sum_n \abs{\widehat{f}(n)}^2 = \int_{\T} \abs{f(\zeta)}^2 dm(\zeta) < \infty,
\]
where $dm$ denotes the unit-normalized arc-length measure on $\T$. We are interested in exploring generic supports and regularity properties of elements in $\ell^2(w)$, as $w$ varies. For instance, if $w_n \downarrow 0$, then one expects $\ell^2(w)$ to contain very singular elements. In fact, if $w_n \downarrow 0$ then the corresponding space $\ell^2(w)$ contains distributions, which are not finite Borel measures on $\T$, see \cite{limani2025generic}. In the opposite regime, when $\lambda=(\lambda_n)_n$ is an increasing sequence, then the corresponding space $\ell^2(\lambda)$ should contain better behaved functions. For instance, if $\sum_n \frac{1}{\lambda_n} < \infty$, then the Cauchy-Schwarz inequality implies that
\[
\sum_n \abs{\widehat{f}(n)} \leq \norm{f}_{\ell^2(\lambda)} \left( \sum_n \frac{1}{\lambda_n} \right)^{1/2}.
\]
In particular, this means that elements in $f\in \ell^2(\lambda)$ are always continuous on $\T$. The following result was recently observed by the author, showing that the above simple observation is essentially sharp.

\begin{thm}[See \cite{limani2025generic}]\thlabel{THM:BADSUPP} Let $(\lambda_n)_n$ be positive numbers satisfying the following conditions:
\begin{enumerate}
    \item[(i)] $\sum_n \frac{1}{\lambda_n} = +\infty$,
    \item[(ii)] there exists $C>1$, such that for all $n\geq 1$:
    \[
    C^{-1} \lambda_n \leq \lambda_k \leq C \lambda_n, \qquad n\leq k \leq 2n.
    \]
\end{enumerate}
Then there exists a positive function $f\in L^\infty(\T,dm)$ such that $\supp{fdm}$ contains no interior, and
\[
\sum_n \abs{\widehat{f}(n)}^2 \lambda_{|n|} < \infty.
\]
\end{thm}
In other words, one can exhibit elements with maximally small topological support, if the sum appearing in $(i)$ diverges. Surprisingly, it was further observed that the additional regularity hypothesis in $(ii)$ cannot be simply removed, nor can it be substituted by a monotonicity condition.

\begin{thm}\thlabel{THM:(ii)} There exists positive numbers $(\lambda_n)_n$ with $\lambda_n \uparrow +\infty$ and $\sum_n \frac{1}{\lambda_n} = + \infty$, such that whenever $S$ is a non-trivial distribution in $\T$ with 
\[
\sum_n \abs{\widehat{S}(n)}^2 \lambda_{|n|}< \infty,
\]
then $\supp{S}= \T$.
\end{thm}
Results of such kind somewhat complement earlier works of T. W. K\"orner on Theorems of Ivashev-Musatov flavor, see \cite{korner1977theorem2} and \cite{korner2003topological}. We also refer the reader to the excellent book of J. P. Kahane in \cite{kahane2019series} on various topics of similar flavor.

\subsection{Purpose and aim}
It was further observed in \cite{limani2025generic} that methods involving sparse Fourier support (Fourier spectrum), such as Sidon sets, Riesz products, and lacunary Hadamard series, cannot produce functions appearing in \thref{THM:BADSUPP}. Our main purpose in this note is to prove thresholds of similar type, which principally rely on classical methods of sparse Fourier support. Even though we confine ourselves within the regime of Hilbert spaces, our methods do not intrinsically rely on this structure in any crucial way, and results of similar flavor have recently also been investigated in \cite{limani2025summable}. Our main results are two-fold and will be phrased in the following section. We shall exhibit thresholds for when one can exhibit elements with small support and with bad range.



\section{Main results}

    




\subsection{Positive measures with small support}
Our main purpose in this section is to determine when elements with arbitrary small support are "typical" in $\ell^2(w)$. To this end, we introduce the complete metric space $(\mathscr{C}, d_{\mathscr{C}})$, where $\mathscr{C}$ is the set of non-empty compact subsets of $\T$, equipped with the so-called \emph{Hausdorff metric}
\[
d_{\mathscr{C}}(E,F) := \sup_{\zeta \in E} \dist{\zeta}{F} + \sup_{\xi \in F} \dist{\xi}{E}.
\]
Let $\mathscr{P}_{w}$ denote the sub-collection of order pairs $(\sigma,E) \in \ell^2(w) \times \mathscr{C}$ satisfying:
\begin{enumerate}
    \item[(a.)] $\supp{\sigma} \subseteq E$,
    \item[(b.)] $\sigma$ is a positive finite Borel measure on $\T$.
\end{enumerate}
Equipping the set $\mathscr{P}_{w}$ with the metric
\[
d_{w}\left( (\sigma,E), (\tau,F) \right) = d_{\mathscr{C}}(E,F) + \norm{\sigma-\tau}_{\ell^2(w)},
\] 
one obtains a complete metric space $(\mathscr{P}_w, d_w)$, see \thref{LEM:FUNCAN} below. Recall that a subset $E$ of a metric space $X$ is said to be \emph{generic} in $X$ if its complement is a meager set, meaning a countable union of nowhere dense sets in $X$. Our first main result asserts that, under natural assumptions on $w$ the set of elements in $\ell^2(w)$ with arbitrarily small support is generic in the following precise sense.

\begin{thm}\thlabel{THM:BAIREPROB} Let $w=(w_n)_{n=0}^\infty$ be positive numbers satisfying the following properties:
\begin{enumerate}
    \item[(i)] $\liminf_{n} w_n = 0$,
    \item[(ii)] $(w_n)_n$ is non-increasing, or there exists $C>1$, such that for any $n\geq 1$:
    \[
    C^{-1} w_n \leq w_k \leq C w_n, \qquad n\leq k\leq 2n.
    \]
\end{enumerate}
Then the sub-collection of $\mathscr{P}_w$ consisting of positive finite Borel measures supported on a sets of Lebesgue measure zero, is generic in the complete metric space $(\mathscr{P}_w, d_w)$.
\end{thm}
Note that the condition $(i)$ and the second statement in $(ii)$ do not necessarily imply that $w_n \to 0$. Our proof of \thref{THM:BAIREPROB} relies on machinery from earlier works of T. W. K\"orner in \cite{korner1986theorem3}, on Ivashev-Musatov type Theorems, and are closer in spirit to the authors recent work in \cite{limani2025generic} and in \cite{limani2025summable}. Under certain natural regularity assumptions on the sequence $w$, the construction of a probability measures with small support can also easily be derived from the abstract potential theoretical framework developed by L. Carleson in his classical survey \cite{carleson1967selected}. See also the work of P. Malliavin in \cite{malliavin1977analyse} on related matters. 

Now the statement of \thref{THM:BAIREPROB} naturally raises the question of whether the hypothesis $(i)$-$(ii)$ are optimal. Notably, the necessity of $(i)$ follows from the observation that if $\liminf_n w_n >0$ then \( \ell^2(w) \hookrightarrow L^2(\mathbb{T}, dm) \), thus every non-zero element in \( \ell^2(w) \) must be supported on a set of positive Lebesgue measure. Next, we demonstrate that the condition $\liminf_n w_n =0$ alone is far from sufficient.

\begin{thm}\thlabel{THM:NECREG} There exists positive numbers $(w_n)_n$ with $\liminf_n w_n =0$, such that any non-trivial complex finite Borel measure $\mu$ on $\T$ with $ (\widehat{\mu}(n))_n \in \ell^2(w)$ has full support.
\end{thm}

Note that the sequence appearing in the statement cannot possibly be bounded, since that would force the containment of $L^2(\T,dm)$, which certainly inhabits elements which do not full support. Despite this negative result, we shall at the end of Section 4 construct a measure with small support in $\ell^2(w)$ under the condition that $(w_n)_n$ tends to zero along a subsequence of long blocks, see \thref{THM:singl2w} below. It will be clear from the proof that the sequence $(w_n)_n$ may be chosen to tend to zero at any prescribed rate and along any lacunary subsequence, provided that $(w_n)_n$ is sufficiently large outside that subsequence. Indeed, our proof of \thref{THM:NECREG} is based on a uncertainty principle, which roughly asserts that the support of a measure cannot be well-localized if its mass is mostly carried by low-order frequencies and a lacunary subset of large-frequencies. As this observation may be of independent interest, we include it here. 
\begin{prop}\thlabel{LEM:UPZYG} For any integer $N\geq 1$ and any numbers $0<\gamma,\delta <1$, there exists a small number $\varepsilon=\varepsilon(N,\gamma,\delta)>0$, such that the following statement holds: whenever $\mu$ is a complex finite Borel measure on $\T$ with total variation norm $\leq 1$ and
\begin{enumerate}
    \item[(i)] $\sum_{|n|\leq N} \abs{\widehat{\mu}(n)}^2 \geq \gamma$,
    \item[(ii)] $\abs{\widehat{\mu}(n)} \leq \varepsilon$, for all $|n|>N$, except for possibly a $3$-sparse set.
\end{enumerate}
then $\sup_{\zeta\in \T} \dist{\zeta}{\supp{\mu}}\leq \delta$.
\end{prop}
Here, subset of integers $\Lambda \subset \mathbb{Z}$ is said to be $N$-sparse for some $N>1$, if 
\[
\inf_{\substack{n,m\in \Lambda \\ |n|>|m|}} \frac{|n|}{|m|}\geq N.
\]
This result may be viewed as a fusion of Lemma 2.1 in \cite{korner2003topological} by T. W. K\"orner (whose proof is attributed to Y. Meyer), and certain uniqueness principle for measure with lacunary Fourier spectrum, which originates from the work of A. Zygmund, see \thref{LEM:ZYGlac} for a precise statement. Before closing this subsection, we make one final observation. One may wonder whether the space $\ell^2(w)$ in \thref{THM:NECREG} perhaps contain a distribution with small support, which would imply that there are $\ell^2(w)$-spaces which exhibit the so-called Piatetski-Shapiro phenomenon: there exist a compact set which supports a distribution in $\ell^2(w)$, but no such measure. This outstanding problem is likely extremely difficult in this general regime, keeping in mind that the classical potential theory framework ensures that this cannot happen for weights of the form $w_n = (1+n)^{-\alpha}$ with $0<\alpha\leq 1$. For recent consideration of problems of this type, we refer the reader to the work of N. Lev and A. Olevskii in \cite{lev2005piatetski} and \cite{lev2011wiener}.



\subsection{Functions with topologically wild range}
In this subsection, we investigate the critical regime for when one can exhibit elements in $\ell^2(\lambda)$ with maximally bad topological range, which complements \thref{THM:BADSUPP} on elements with maximally small topological support. To this end, we recall that a measurable function $f$ on $\T$ is said to be \emph{everywhere locally unbounded}, if for any arc $I \subset \T$ and any number $M>0$, we have 
\[
m \left( \{\zeta \in I: \abs{f(\zeta)} >M \} \right) >0.
\]
To introduce the appropriate setting, we consider the modified Hilbert space $\mathscr{S}^2(\lambda)$ of functions $f\in \ell^2(\lambda)\cap L^2(\T,dm)$ normed by
\[
\norm{f}_{\mathscr{S}^2(\lambda)} := \norm{f}_{\ell^2(\lambda)} + \norm{f}_{L^2}.
\]
Our main result in this subsection is the following complete characterization of functions with maximally wild range. 

\begin{thm}\thlabel{THM:BAIRETOP} Let $(\lambda_n)_n$ be a sequence of positive numbers. Then the following statements are all equivalent:
\begin{enumerate}
    \item[(i)] $\mathscr{S}^2(\lambda)$ contains an everywhere locally unbounded function,
    \item[(ii)] the collection of everywhere locally unbounded functions in $\mathscr{S}^2(\lambda)$ is generic, 
    \item[(iii)] $\sum_n \frac{1}{\lambda_n} = + \infty$.
\end{enumerate}
\end{thm}
A couple of remarks regarding \thref{THM:BAIRETOP} are now in order. First, the reason for introducing the space $\mathscr{S}^2(\lambda)$ is to avoid the degenerate behavior $\liminf_n \lambda_n =0$, hence ensure that $\ell^2(\lambda)$ only consists of $L^2(\T,dm)$-functions, for which we always have that $|f|< \infty$ $dm$-a.e on $\T$. Actually, that degenerate case can be handled separately (see Step 1 in the proof of \thref{LEM:l2C}, by means of constructing a lacunary Hadamard series which is locally unbounded everywhere. The second remark concerns the proof of \thref{THM:BAIRETOP}, which in this full generality is non-constructive, and relies on certain duality arguments. However, it turns out that for a fairly broad class of weight sequences $(\lambda_n)_n$, it is actually possible to give a purely constructive proof. 
\begin{thm}\thlabel{THM:BADRANGE} Let $(\lambda_n)_n$ be positive numbers satisfying the following hypothesis:
\begin{enumerate}
    \item[(i)] $\sum_n \frac{1}{\lambda_n} = +\infty$,
    \item[(ii)] $(\lambda_n)_n$ is non-decreasing, or there exists $C>1$, such that for all $n\geq 1$:
    \[
    C^{-1} \lambda_n \leq \lambda_k \leq C \lambda_n, \qquad n\leq k \leq 2n.
    \]
\end{enumerate}
Then a positive $f\in L^2(\T,dm)$ can be constructed, such that $f$ is everywhere locally unbounded and satisfies
\[
\sum_n \abs{\widehat{f}(n)}^2 \lambda_{|n|} < \infty.
\]
\end{thm}
The proof of \thref{THM:BADRANGE} involves constructing a special family of Riesz products, which will be outlined in the first part of Section 3. Notably, \thref{THM:BADRANGE} holds for the class of fractional Sobolev spaces consisting of functions \( f \in L^2(\mathbb{T}, dm) \) such that
\[
\sum_{n \in \mathbb{Z}} |\widehat{f}(n)|^2 (1 + |n|)^s < \infty, \qquad 0 < s \leq 1.
\]
However, in this classical regime, such observations are certainly not new, and well-known by experts. For instance, in the case when \( \lambda_n = 1 + n \), any function $f\in \ell^2(\lambda)$ must be finite outside a set of logarithmic capacity zero, and one can construct such a function $f$ which blows-up on every point of a fixed compact set of logarithmic capacity zero. Similar conclusions also hold for weights \( \lambda_n = (1 + n)^s \) with \( 0 < s < 1 \), where the exceptional sets are governed by Riesz capacities. We refer the reader to the survey of L. Carleson in \cite{carleson1967selected} for details on such classical results, and further generalizations beyond.

\subsection{Organization and notation.}
We conclude this section with a brief overview of the manuscript. Section 3 is devoted to proving our main results announced in Subsection 2.2, with principal emphasis on constructing elements appearing in the statement of \thref{THM:BADRANGE}. The last subsection contains \thref{THM:BAIRETOP} on generic elements in weighted $\ell^2$ with maximally bad topological range. The proofs of our results on elements with small support announced in Subsection 2.1, are contained in Section 4. The subsections 4.1-4.3 are devoted to \thref{THM:BAIREPROB}, while the proof of \thref{THM:NECREG}, based on \thref{LEM:UPZYG}, are both located in subsection 4.4. At last, we complete Section 4 with a short and constructive proof of \thref{THM:singl2w}. 

Throughout, we shall frequently utilize the notation $A \lesssim B$ to designate that two positive numbers $A, B >0$ satisfy the relation $A \leq cB$, for some positive constant $c>0$. In case $A\lesssim B$ and $B \lesssim A$ hold, we shall often write $A\asymp B$. Additionally, absolute constants will generally be denoted by $C,c$, and alike, but their exact value may vary from line to line. In case the dependence of a constant $C$ on specific parameter(s), say $\varepsilon$, needs to be stressed, we shall designate it by $C(\varepsilon)$.

\subsection{Acknowledgments} This research was supported by the generous stipend, provided by the Knut \& Alice Wallenberg Foundation (grant no. 2021.0294). The author would like to thank Andrei Semenov for useful feedback on an earlier version of this manuscript.

\section{Functions with wild range}\label{SEC:3}
Here, we shall illustrate how techniques involving spectral gaps, such as a family of Riesz products, can be utilized in order to prove \thref{THM:BADRANGE}. 

\subsection{Riesz products}
Let $(N_j)_j$ be a lacunary sequence of positive integers satisfying $\kappa := \inf_j N_{j+1}/N_j \geq 3$, and let $(a_j)_j$ be positive real numbers taking values between $0<a_j<1$. The classical Riesz partial products are defined as
\[
d\sigma_N(\zeta) = \prod_{j=1}^N \left( 1+ a_j \Re(\zeta^{N_j}) \right) dm(\zeta), \qquad \zeta \in \T, \qquad N=1,2,3,\dots
\]
It is well-known that the assumptions on $(a_j)_j$ ensure that $\sigma_N$ converges to a probability measure $\sigma$ on $\T$, in the weak-star topology of finite Borel measure on $\T$. It is therefore convenient to refer to the limiting measure $\sigma$ as the Riesz product associated to the parameters $(a_j)_j$ and $(N_j)_j$. Now a crucial property of Riesz products is that the lacunary assumption on $(N_j)_j$ with $\kappa \geq 3$ ensures that Fourier support of $\sigma$ is contained in the union of disjoint blocks of the form 
\begin{equation}\label{EQ:Lambdaj}
\Lambda_j := \left\{ (1-(\kappa-1)^{-1}) N_j \leq |n| \leq (1+(\kappa-1)^{-1}) N_j \right\}, \qquad  j=0,1,2,\dots,
\end{equation}
where $\Lambda_0 := \{1\}$ by default. Using this observation, it is straightforward to deduce the Fourier estimates
\[
\sum_{n\in \Lambda_{j}} \abs{\widehat{\sigma}(n)}^2 \leq a^2_{j} \prod_{k=1}^j (1+a^2_k), \qquad j=1,2,3,\dots.
\]
In fact, one can show that $\sum_j a^2_j < \infty$ if and only if $d\sigma = f dm$ with $f\in L^2(\T, dm)$. For further details on this matter, we refer the reader Ch. IV in \cite{katznelson2004introduction}. A deeper property that will be utilized in the next section to prove \thref{THM:singl2w}, is that if $\sum_j a^2_j = \infty$ and $a_j \to 0$, then the corresponding Riesz product $\sigma$ must be supported on a set of Lebesgue measure zero in $\T$, see Theorem 7.7 from Ch. V in \cite{zygmundtrigseries}.

\subsection{An adapted family of Riesz products}
For the purpose of establishing \thref{THM:BADRANGE}, we shall actually need to introduce a class of Riesz products, specifically adapted to our framework. To this end, we consider the normalized Dirichlet-type kernels $(T_N)_N$ defined as 
\begin{equation}\label{EQ:DIRI}
T_N(\zeta) = \frac{1}{N} \sum_{0\leq k< N} \zeta^{N+k}, \qquad \zeta \in \T, \qquad N=1,2,3,\dots
\end{equation}
Below, we record a couple of standard properties of $(T_N)_N$, which are straightforward to verify:
\begin{enumerate}
    \item[(a.)] $\int_{\T} T_N dm =0$,
    \item[(b.)] $\supp{\widehat{T_N}} \subseteq \{N\leq |n|< 2N\}$,
    \item[(c.)] $\sup_{\zeta \in \T} \abs{T_N(\zeta)} = T_N(1)=1$,
    \item[(d.)] $\norm{T_N}^2_{L^2} = 1/N$.
\end{enumerate}
Using these normalized Dirichlet-type kernels as our principal building-blocks, we form the corresponding Riesz-type products as follows:
\begin{equation} \label{EQ:RIESZPART}
P_n(\zeta) = \prod_{j=1}^n \left( 1 + a_j \Re(T_{N_j}(\zeta) \right), \qquad \zeta \in \T, \qquad n=1,2,3,\dots,
\end{equation}
where $(a_j)_j$ are positive numbers with $0<a_j \leq 1/2$, and $(N_j)_j$ positive lacunary integers with $\inf_j N_{j+1}/N_j \geq 4$, say. Since each $T_N$ have zero average on $\T$, as demonstrated in $(a.)$, and $0<a_j<1$, it is easy to see that $d\nu_n := P_n dm$ are probability measures for each $n$. It follows from a standard weak-compactness argument involving functional analysis, that $\nu_n$ converge (perhaps up to an appropriate subsequence) to a probability measure $\nu$ on $\T$, which we shall refer to as the \emph{Riesz product} adapted to $(T_{N_j})_j$ and $(a_j)_j$. Our first lemma in this subsection contains the important features about the Fourier coefficients of these adapted Riesz products, that will be utilized in the next subsection.

\begin{lemma}\thlabel{LEM:RIESZPROD} Let $(a_j)_j$ be positive numbers with $0<a_j \leq 1/2$ and let $(N_j)_j$ be positive integers with $\inf_j N_{j+1}/N_j \geq 4$. Consider the Riesz product adapted to $(T_{N_j})_j$ and $(a_j)_j$, defined as 
\[
d\nu(\zeta) := \prod_{j=1}^\infty \left( 1 + a_j \Re(T_{N_j}(\zeta) \right)dm(\zeta), \qquad \zeta \in \T,
\]
interpreted in sense of weak-star limits of finite Borel measures on $\T$. Then the following statements hold: 
\begin{enumerate}
    \item[(i)] Then $d\nu = fdm$ with a positive density $f$ belonging to $L^2(\T,m)$.
    \item[(ii)] The Fourier support of $\nu$ is contained in union of sparse sets $\Lambda_j := \{ 3^{-1} N_j \leq |n| \leq 3N_j \}$, $j=0,1,2,\dots$, where $\Lambda_0:=\{1\}$ by default. Moreover, the following local $\ell^2$-estimates hold:
    \[
    \sum_{n\in \Lambda_{j+1}} \abs{\widehat{\nu}(n)}^2 \leq \frac{a^2_{j+1}}{N_{j+1}} \prod_{k=1}^j \left(1 + \frac{a^2_k}{N_k} \right), \qquad j=1,2,3,\dots
    \]

\end{enumerate}

\end{lemma}

\begin{proof} 

Note that the corresponding Riesz partial products in \eqref{EQ:RIESZPART} satisfy the property
\begin{equation}\label{EQ:Recur}
P_{j+1}(\zeta) = P_j(\zeta) + a_{j+1} P_j(\zeta) \Re T_{N_{j+1}}(\zeta), \qquad j=1,2,3, \dots.
\end{equation}
Note that the property $(b.)$ of the Dirichlet-type kernels ensure that the Fourier spectrum of $P_1 := 1 +a_1 \Re(T_{N_1})$ is contained in the set of integers $\{1\} \cup \{N_1\leq |n|\leq 2N_1\}$. By means of induction, suppose that the Fourier support of $P_j$ is contained in the union $\cup_{k=0}^j \Lambda_k$. Now in views of \eqref{EQ:Recur}, it thus suffices to verify that $R_{j+1} := P_j \Re T_{N_{j+1}}$ has Fourier support in $\Lambda_{j+1}$. Observe that a successive application of the lacunary assumption $\inf_j N_{j+1}/N_j \geq 4$, implies that $N_{j+1} \geq 3(N_j + \dots + N_1)$ for all $j$, hence a simple calculation gives 
\[
\supp{\widehat{R_{j+1}}} \subseteq \{ N_{j+1} - \sum_{k=1}^j N_k \leq |n| \leq 2N_{j+1} + \sum_{k=1}^j N_k \} \subseteq \Lambda_{j+1}.
\]
This proves that each Riesz partial product $P_j$ has Fourier support in $\cup_{k=0}^j \Lambda_k$, hence the claim for the support of $\sigma$ readily follows. To prove the local $\ell^2$-estimates, we primarily observe that 
\[
\sum_{n\in \Lambda_1} \abs{\widehat{\sigma}(n)}^2 = \sum_{n\in \Lambda_1} \abs{\widehat{P_1}(n)}^2 = a^2_1\norm{T_{N_1}}^2_{L^2} = \frac{a^2_1}{N_1}.
\]
It follows from \eqref{EQ:Recur} and the Fourier supports of $P_1$ and $\Re(T_{N_2})$, that
\[
\sum_{n\in \Lambda_2} \abs{\widehat{\sigma}(n)}^2 = \sum_{n\in \Lambda_2} a^2_2 \abs{\widehat{P_1 \cdot \Re T_{N_2}}(n)}^2 \leq a^2_2 \norm{P_1}^2_{L^2} \norm{T_{N_2}}^2_{L^2} = \frac{a^2_2}{N_2} (1 + \frac{a^2_1}{N_1}). 
\]
Now a straightforward induction argument allows us to upgrade this estimate for all $j\geq 1$:
\[
\sum_{n\in \Lambda_{j+1}} \abs{\widehat{\sigma}(n)}^2 \leq \frac{a^2_{j+1}}{N_{j+1}} \prod_{k=1}^j \left(1 + \frac{a^2_k}{N_k} \right), 
\]
thus establishing the claim in $(ii)$. 

In order to prove $(i)$, we note that the by means of summing up the local $\ell^2$-estimates in $(ii)$ and using \eqref{EQ:Recur}, we see that $\sigma \in L^2$ if the infinite product $\prod_{j}(1+a^2_j/N_j)$ converges. But this is indeed convergent in view of $0<a_j<1$ and $\inf_j  N_{j+1}/N_j \geq 4$.

\end{proof}

Our next lemma illustrates an important uniqueness principle for the adapted family of Riesz products introduced in \eqref{EQ:RIESZPART}. The corresponding results are well-known for classical Riesz products (for instance, see Ch. VI, Section 6, Theorem 6.3 in \cite{zygmundtrigseries}, and Ch. V, Exercise for Section 1 in \cite{katznelson2004introduction}). While our proof is inspired from classical arguments, it requires several subtle modifications tailored to the present setting.

\begin{lemma}\thlabel{LEM:RIGIDLAC} Let $(a_j)_j$ be positive numbers with $0<a_j\leq 1/2$ and $\sum_j a^2_j <\infty$, and let $(N_j)_j$ be positive integers satisfying $\inf_j N_{j+1}/N_j \geq 4$. If the Riesz product $P$ adapted to $(a_j)_j$ and $(N_j)_j$ is essentially bounded on an arc $I\subseteq \T$, then $\sum_j a_j < \infty$.
\end{lemma}
\begin{proof} Note that according to $(i)$ of \thref{LEM:RIESZPROD}, we have that $P$ is a positive function in $L^2(\T, dm)$, therefore we can write
\[
\log P(\zeta) = \sum_{j=1}^\infty \log \left( 1+ a_j \Re T_{N_j}(\zeta) \right), \qquad \zeta \in \T.
\]
Indeed, the simple estimate $\abs{\log(1+x)}\leq |x|$ for $|x|\leq 1/2$ in conjunction with the assumptions $0<a_j<1/2$, $\sum_j a^2_j < \infty$, and $\norm{T_N}_{L^\infty}\leq 1$ implies that the above identity holds in the sense of convergence in $L^2(\T,dm)$. In fact, the improved estimate $\abs{\log(1+x)-x}\leq x^2/2$ for $|x|\leq 1/2$, in conjunction with $\sum_j a^2_j$ actually implies that we can write
\[
\log P(\zeta) + h(\zeta) =  \sum_{j=1}^\infty a_j \Re T_{N_j}(\zeta), \qquad \zeta \in \T
\]
for some bounded real-valued function $h$ on $\T$. Arguing by contraposition, we assume that there exists a (closed) arc $I_0 \subseteq \T$, such that $P$ is essentially bounded on $I_0$. Setting $F:= \log P + h$, this implies that $F\in L^2(\T,dm)$ and is essentially bounded above on $I_0$, that is, there exists a number $M_0>0$, such that 
\[
m(I_0 \cap \{F>M_0 \} ) =0.
\]
To this end, consider the classical family of finite Riesz products
\[
Q_n(\zeta) := \prod_{j=1}^n \left( 1 + \Re \left( \zeta^{N_j} \right) \right), \qquad \zeta \in \T, \qquad n=1,2,3,\dots
\]
which are positive and satisfy $\int_{\T} Q_n dm =1$. Pick a non-negative smooth function $\phi_0$, which is supported on $I_0$ and satisfies $\widehat{\phi_0}(0)=1$, and note that since $\phi_0 Q_n$ is non-negative on $\T$, we have
\begin{equation}\label{EQ:F_+}
\int_{\T} \phi_0(\zeta) Q_n(\zeta) F(\zeta) dm(\zeta) \leq \norm{F_+}_{L^\infty(I_0)} \int_{I_0} \phi_0 Q_n dm \leq \norm{F_+}_{L^\infty(I_0)} \norm{\phi_0}_{L^\infty(I_0)},
\end{equation}
where $F_+ := \max (F, 0)$. On the other hand, a straightforward expansion implies that 
\[
\int_{\T} \phi_0(\zeta) Q_n(\zeta) F(\zeta) dm(\zeta) = \sum_{j=1}^\infty a_j \Re \int_{\T} \phi_0 Q_n T_{N_j} dm = \sum_{j=1}^\infty a_j \frac{1}{N_j} \sum_{1\leq k\leq N_j} \widehat{\phi_0 Q_n}(k+N_j).
\]
Now if terms of the sum on the right hand side would consist of the Fourier coefficients of $Q_n$ instead of  $\phi_0 Q_n$, then we would easily compute and estimate the sum as follows:
\begin{equation}\label{EQ:ckappa}
\sum_{j=1}^\infty a_j \frac{1}{N_j} \sum_{1\leq k\leq N_j} \widehat{Q_n}(k+N_j) \geq c(\kappa) \sum_{j=1}^{n} a_j,
\end{equation}
where we used the standard fact that there exists a constant $0<c(\kappa)<1$, only depending on $\kappa:= \inf_j N_{j+1}/N_j$, such that $\widehat{Q_n}(k+N_j) =1$ for $0\leq k \leq c(\kappa)N_j$, and zero else. Therefore, it only remains to show that, say 
\[
\abs{\widehat{\phi_0 Q_n}(k) - \widehat{Q_n}(k)} \leq \frac{c(\kappa)}{10}, \qquad N_{k_0} \leq k \leq N_{n+1}, 
\]
where the integer $k_0>0$, may depend on $\phi_0$ and $c(\kappa)$. However, using the fact that $\widehat{\phi_0}(0)=1$, this is straightforward to verify as follows:
\begin{multline*}
\abs{\widehat{\phi_0 Q_n}(k) - \widehat{Q_n}(k)} \leq \sum_{j\neq k} \abs{\widehat{\phi_0}(k-j)} \abs{\widehat{Q_n}(j)}  \\ \lesssim \norm{\phi''_0}_{L^\infty} \sum_{\substack{j\neq k \\ j\in \supp{\widehat{Q_n}}}} |j-k|^{-2} 
\lesssim \norm{\phi''_0}_{L^\infty} \sum_{|j|\geq N_{k_0}} \frac{1}{j^2} \leq \frac{c(\kappa)}{10},
\end{multline*}
whenever $k \geq N_{k_0}$ with $k_0>0$ large enough. Now this in conjunction with \eqref{EQ:ckappa} implies that 
\[
\sum_{j=k_0}^\infty a_j \frac{1}{N_j} \sum_{1\leq k\leq N_j} \widehat{\phi_0 Q_n}(k+N_j) 
\geq \frac{9c(\kappa)}{10} \sum_{k_0\leq j\leq n} a_j.
\]
Returning back to \eqref{EQ:F_+} and letting $n\to \infty$, we arrive at 
\[
\sum_{j>k_0} a_j \leq C(I_0, \kappa) \norm{F_+}_{L^\infty(I_0)}
\]
where $C(I_0, \kappa)>0$ is a constant only depending on $I_0$ and $\kappa$. This completes the proof.
\end{proof}

We remark that the above proof actually implies that $\sum_j a_j < \infty$ if and only if $F = \sum_j a_j \Re T_{N_j}$ is essentially bounded above on some/any arc $I\subset \T$. In particular, if $F$ vanishes on an arc, then it vanishes identically on $\T$.

\subsection{The construction of elements with maximally bad range}
Using the developments from the previous subsection, we shall here give a constructive proof of \thref{THM:BADRANGE}. Before turning our attention to this primary task, we need one additional lemma on summable sequences, whose simple proof is sketched below.

\begin{lemma}\thlabel{LEM:SEQ1}  Let $(\lambda_n)_n$ be positive numbers satisfying the hypothesis $(i)-(ii)$ of \thref{THM:BADRANGE}. There exists positive numbers $(a_n)_n$ with $0<a_n <1/2$ that satisfy the following properties:
\[
(1.) \, \sum_{n>0} a_n = + \infty, \qquad (2.) \, \sum_{n>0} \max \left(1, \lambda_{4^n} 4^{-n} \right) a^2_{n} < \infty.
\]
\end{lemma}

\begin{proof} We primarily observe that the hypothesis $(ii)$ on $(\lambda_n)_n$ allows one to rephrase $(i)$ equivalently as 
\begin{equation}\label{EQ:CAUCHYCOND}
\sum_j \frac{4^j}{\lambda_{4^j}} \asymp \sum_{j\geq 0} \sum_{4^j\leq |n| < 4^{j+1}}\frac{1}{\lambda_n} = \sum_n \frac{1}{\lambda_n} =+ \infty.
\end{equation}
With this observation in mind, we set $d_j := 4^j \lambda^{-1}_{4^j}$ and note that by means of substituting $d_j$ with $\min(1, d_j)$, we may without loss of generality assume that $d_j \leq 1$. Therefore the desired statement reduces to proving:
\[
(1.) \, \, \sum_j a_j <\infty, \qquad (2.) \, \, \sum_j \frac{a^2_j}{d_j} < \infty.
\]
To this end, consider the sequence of real numbers defined by
\[
a_j = \frac{1}{2}\frac{d_j}{1+ \sum_{k=1}^j d_k}, \qquad j=1,2,3, \dots,
\]
which readily satisfy $0<a_j\leq 1/2$. Since $\sum_j d_j = +\infty$, by assumption, a straightforward argument involving the integral test also implies that 
\[
\sum_j a_j = \frac{1}{2} \sum_j \frac{d_j}{1+ \sum_{k=1}^j d_k} =+ \infty.
\]
On the other hand, a similar argument also shows that
\[
\sum_j \frac{a^2_j}{d_j} =\frac{1}{4} \sum_j \frac{d_j}{(1+ \sum_{k=1}^j d_k)^2} < \infty.
\]

\end{proof}

\begin{proof}[Proof of \thref{THM:BADRANGE}] 
Let $(a_j)_j$ be a collection of positive real numbers appearing in the statement of \thref{LEM:SEQ1}, and consider the associated Riesz product, adapted to $(T_{4^j})_j$ and $(a_j)_j$, defined as
\[
d\nu(\zeta) := \lim_{N\to \infty} \prod_{j=1}^N \left( 1 + a_j \Re( T_{4^j}(\zeta) \right) dm(\zeta), \qquad \zeta \in \T,
\]
interpreted in sense of weak-star limit of finite measure on $\T$. Now the property $(2.)$ of $(a_j)_j$ from \thref{LEM:SEQ1} in conjunction with \thref{LEM:RIESZPROD} $(a.)$ implies that $d\nu = fdm$, where the non-negative density $f$ belongs to $L^2(\T, dm)$. In fact, the using the local $\ell^2$-estimates in part $(b.)$ of \thref{LEM:RIESZPROD}, we get 
\[
\sum_{n\in \Lambda_{j+1}}  \abs{\widehat{f}(n)}^2 \lesssim \frac{a^2_{j+1}}{4^{j+1}} \prod_{1\leq k\leq j} \left( 1 +\frac{a^2_{k}}{4^k} \right) \asymp \frac{a^2_{j+1}}{4^{j+1}}, \qquad j=0,1,2 \dots
\]
where the last inequality follows from the observation that $0<a_j<1$ and that the infinity product $\prod_{j} (1+ a^2_j 4^{-j})$ converges. With this estimate at hand, and in view of $(2.)$ from \thref{LEM:SEQ1}, we get
\[
\sum_{n\neq 0} \abs{\widehat{f}(n)}^2 \lambda_{|n|} \asymp \sum_{j\geq 0} \lambda_{4^j} \sum_{n\in \Lambda_{j+1}}  \abs{\widehat{f}(n)}^2 \asymp \sum_j \lambda_{4^j} \frac{a^2_j}{4^j} < \infty,
\]
hence $(\widehat{f}(n))_n$ belongs to $\ell^2(\lambda)$. Meanwhile, the property $(1.)$ of $(a_j)_j$ from \thref{LEM:SEQ1} in conjunction with \thref{LEM:RIGIDLAC} also implies that $f$ is essentially unbounded in every neighborhood of $\T$. This completes the proof. 
\end{proof}

We can now explain the necessity of introducing the family of adapted Riesz products. If one attempts to construct an element $f$ as in the statement of \thref{THM:BADRANGE} using classical families of Riesz products, then the properties $(\widehat{f}(n))_n \in \ell^2(\lambda)$ and $f$ being locally everywhere unbounded, are seen to be equivalent to the requirements
\[
(a.) \qquad \sum_n a_n^2 \lambda_{4^n} < \infty, \qquad (b.) \qquad \sum_n a_n = \infty.
\]
However, if $\lambda_n = 1+n$, for instance, a straightforward application of the Cauchy–Schwarz inequality shows that these requirements are mutually incompatible. The main reason for introducing the adapted Riesz products via normalized Dirichlet-type kernels is due to their accelerated $L^2$-decay, which transform $(a.)$ into the condition $\sum_n a_n^2 \lambda_{4^n}4^{-n} < \infty$, enabling us to reach the critical threshold.

\subsection{Functional analytic setting and duality}

In this subsection, we demonstrate how to dispense with the additional regularity condition $(ii)$, and thus prove \thref{THM:BAIRETOP} in full generality. Unfortunately, the our proof will no longer be constructive, as it will be accomplished using functional analysis. We shall start by setting up the basic functional analytic framework that will be relevant for our further developments. To this end, we recall that for any sequence of positive real numbers $(\lambda_n)_n$, the corresponding space $\ell^2(\lambda)$ becomes a Hilbert space wrt to the inner product
\[
\langle f, g \rangle_{\ell^2(\lambda)} := \sum_n \widehat{f}(n) \conj{\widehat{g}(n)} \lambda_{|n|}.
\]
It is easy to see that $\ell^2(\lambda)$ is separable with the trigonometric polynomials forming a dense subset. It will be more convenient to identify the dual space $\ell^2(\lambda)^*$ of continuous linear functionals on $\ell^2(\lambda)$ using a different pairing than the canonical inner product of $\ell^2(\lambda)$. To this end, we denote by $\ell^2(1/\lambda)$ the space of distributions $S$ on $\T$ with 
\[
\norm{S}^2_{\ell^2(1/\lambda)} := \sum_n \abs{\widehat{S}(n)}^2 \frac{1}{\lambda_{|n|}} < \infty.
\]
Cauchy-Schwarz inequality implies that every distribution $S$ with Fourier coefficients in $\ell^2(1/\lambda)$ induces a bounded linear functional $l_S$ on $\ell^2(\lambda)$, considered in the classical $\ell^2$-pairing:
\[
l_S(f) = \sum_n \widehat{f}(n) \conj{\widehat{S}(n)}, \qquad  (\widehat{f}(n))_n \in \ell^2(\lambda).
\]
Conversely, if $l$ is a bounded linear functional on $\ell^2(\lambda)$, then the Riesz representation theorem ensures that there exists a unique element $g \in \ell^2(\lambda)$, such that 
\[
l(f) = \langle f, g \rangle_{\ell^2(\lambda)}, \qquad (\widehat{f}(n))_n \in \ell^2(\lambda).
\]
Now for any trigonometric polynomial $T$, we have 
\[
l(T) = \sum_n \widehat{T}(n) l(\zeta^n) = \sum_n \widehat{T}(n) \conj{\widehat{g}(n)} \lambda_n.
\]
Since the trigonometric polynomials are dense, it follows that the distribution $G$ on $\T$ with $\widehat{G}(n) =\conj{\widehat{g}(n)} \lambda_n$ belongs to $\ell^2(1/\lambda)$ with $\norm{G}_{\ell^2(1/\lambda)}= \norm{g}_{\ell^2(\lambda)}$. We conclude that every bounded linear functional on $\ell^2(\lambda)$ arises from summation against an element $\ell^2(1/\lambda)$ in the classical $\ell^2$-pairing, and $\ell^2(\lambda)^*$ is isometrically isomorphic to $\ell^2(1/\lambda)$. Switching the roles of $(\lambda_n)_n$ and $(1/\lambda_n)_n$, we also conclude that $\ell^2(\lambda)$ is reflexive. With these basic notions at hand, we now turn to the basic building-block in our proof of \thref{THM:BAIRETOP}

\begin{lemma}\thlabel{LEM:l2C} Let $(\lambda_n)_n$ be positive numbers with $\sum_n \frac{1}{\lambda_n} =+ \infty$. For any arc $I \subset \T$, there exists an element $f_I\in L^2(\T,dm)$ with $(\widehat{f_I}(n))_n \in \ell^2(\lambda)$, which is essentially unbounded on $I$:
\[
m\left( \{\zeta\in I: \abs{f_I(\zeta)} >M \} \right) >0, \qquad M>0.
\]
\end{lemma}
\begin{proof}
The proof will be divided into the following three steps:
\proofpart{1}{A simple reduction:}
We shall primarily reduce the analysis to the case $\inf_{n} \lambda_n \geq 1$. Suppose there exists an infinite subsequence $(\lambda_{n_j})_j$ with the property that $0<\lambda_{n_j}<1$ for all $j$. By means of passing to yet another subsequence, we may also assume that $\inf_j n_{j+1}/n_j\geq 3$. In that case, we may construct the desired function $f$ fairly explicitly, using lacunary Hadamard series. To this end, pick a sequence of positive numbers $(a_j)_j$ with the properties
\[
(a.) \, \sum_j a^2_j < \infty, \qquad (b.) \, \sum_j a_j = +\infty,
\]
and consider the associated lacunary Hadamard series:
\[
f(\zeta) := \sum_j a_j \zeta^{\pm n_j}, \qquad \zeta \in \T.
\]
The assumption $(a.)$ ensures that $f \in L^2(\T, dm)$. In fact, since $\lambda_{n_j}<1$ for all $j$, we also have 
\[
\sum_n \abs{\widehat{f}(n)}^2 \lambda_{|n|} \leq 2 \sum_j a^2_j \lambda_{n_j} \leq 2\sum_j a^2_j < \infty
\]
thus $(\widehat{f}(n) )_n \in \ell^2(\lambda)$. On the other hand, the assumption $(b.)$ implies that $f$ is essentially unbounded in every neighborhood of $\T$, according to a classical Theorem of A. Zygmund. For instance, see Theorem 6.1 in Ch. VI of \cite{zygmundtrigseries}, and \thref{LEM:RIGIDLAC} in subsection 3.2. Since a desirable function has been constructed in the case when $\liminf_n \lambda_n <1$, we may now assume that $\lambda_n \geq 1$ for all $n$, hence $\ell^2(\lambda) \subseteq L^2(\T,dm)$.

\proofpart{2}{Existence of unbounded functions:}  We shall first show that $\ell^2(\lambda)$ contains an essentially unbounded function. For for the sake of obtaining a contradiction, assume that $\ell^2(\lambda) \hookrightarrow L^\infty(\T,dm)$, hence by the closed graph Theorem, there exists a constant $C(\lambda)>0$, such that 
\[
\norm{f}_{L^\infty} \leq C(\lambda) \norm{f}_{\ell^2(\lambda)} , \qquad \forall f \in \ell^2(\lambda).
\]
Fix an arbitrary element $g \in L^1(\T,dm)$ and note that for any trigonometric polynomial $T$, we have
\[
\abs{\sum_n \widehat{g}(n) \conj{\widehat{T}(n)} } = \abs{\int_{\T} g(\zeta) \conj{T(\zeta)} dm(\zeta)} \leq \norm{g}_{L^1(\T)} \norm{T}_{L^\infty} \leq C(\lambda) \norm{g}_{L^1(\T)} \norm{T}_{\ell^2(\lambda)}.
\]
From this argument, we see that any element $g\in L^1(\T, dm)$ induces a (densely defined) bounded linear functional on $\ell^2(\lambda)$, hence we conclude that $L^1(\T, dm) \hookrightarrow \ell^2(1/\lambda)$, with 
\[
\sum_n \abs{\widehat{g}(n)}^2 \frac{1}{\lambda_{|n|}} \leq C(\lambda)^2 \norm{g}^2_{L^1}, \qquad \forall g \in L^1(\T,dm).
\]
Now let $J_N$ be arcs centered at $\zeta=1$ of length $1/N$ and consider the simplest family of box-kernels $\varphi_N(\zeta) := N 1_{J_N}(\zeta)$, whose Fourier coefficients are easily computed to be the discrete sinc-functions:
\[
\widehat{\varphi_N}(n) = \frac{\sin(n/2N)}{n/2N} \qquad n\neq 0.
\]
A simple computations shows that there exists a numerical constant $c_0>0$, such $\abs{\widehat{\varphi_N}(n)} \geq 1/2$ for all $1\leq |n|\leq c_0 N$. Using this, we get 
\[
\frac{1}{4} \sum_{1\leq |n|\leq c_0 N} \frac{1}{\lambda_{|n|}} \leq \sum_{n} \abs{\widehat{\varphi_N}(n)}^2 \frac{1}{\lambda_{|n|}} \leq C(\lambda)^2 \norm{\varphi_N}^2_{L^1} = C(\lambda)^2, \qquad N=1,2,3,\dots
\]
However, this contradicts the assumption $\sum_n \frac{1}{\lambda_n}$ diverges. We conclude that there exists an $L^2(\T,dm)$-function $f$ with $(\widehat{f}(n))_n \in \ell^2(\lambda)$, which is essentially unbounded on $\T$.

\proofpart{2}{Unbounded at any prescribed arc:} Using a simple argument involving the pigeon hole principle, we show that our candidate $f$ can be modified to be unbounded on desirable arc $I \subseteq \T$. Fix an arbitrary arc $I \subset \T$, and pick an integer $N>1$ such that $m(I)>1/N$. Now consider a decomposition $I_1,\dots, I_N$ of $\T$ consisting of $N$ arcs, each of normalized length $1/N$, and whose interior are mutually disjoint, and write
\[
m(\{\zeta\in \T: \abs{f(\zeta)} >M \}) = \sum_{j=1}^N m(I_j \cap \{ |f|>M \}).
\]
Since $f$ was proved to be essentially unbounded on $\T$, it follows that the left hand side is positive for any positive integers $M=1,2,3, \dots$. By the pigeonhole principle, there must exist an arc $I_j$, such that $m(I_j \cap \{ |f|>M \})>0$ for infinitely many integers $M$. But since the function
\[
t \mapsto m(I_j \cap \{ |f|>t \}) 
\]
is non-increasing, it actually follows that $m(I_j \cap \{ |f|>M \})>0$ for all $M>0$. Now if $\xi_{I}$ denotes the center of an arc $I$, we readily see that the translate $f_I(\zeta) = f(\zeta \conj{\xi_{I_j} \xi_I}) \in \ell^2(\lambda)$ and is essentially unbounded in $I$. This completes the proof.
 
\end{proof}



\subsection{A Baire category argument}
This subsection contains the proof of \thref{THM:BAIRETOP}. Recall that if $\lambda= (\lambda_n)_n$ is a sequence of positive numbers, we denote by $\mathscr{S}^2(\lambda)$ the Hilbert space of functions $f\in \ell^2(\lambda)\cap L^2(\T,dm)$, equipped with the norm 
\[
\norm{f}_{\mathscr{S}^2(\lambda)} := \norm{f}_{\ell^2(\lambda)} + \norm{f}_{L^2}.
\]
We shall now exhibit a family of meager subsets of $\mathscr{S}^2(\lambda)$, relevant for our purposes. To this end, fix an arc $I \subseteq \T$, a number $M>0$, and consider subsets of $\mathscr{S}^2(\lambda)$ of the form
\begin{equation*}\label{EQ:MEAGER}
\mathscr{L}(I, M) := \left\{ f\in \mathscr{S}^2(\lambda):  m\left( \{\zeta \in I: f(\zeta) \geq M \}\right) = 0\right\}.
\end{equation*}
Our main proposition in thus subsection is essentially contained in the following proposition.

\begin{prop}\thlabel{PROP:BAIRE} Let $(\lambda_n)_n$ be positive numbers satisfying $\sum_n \frac{1}{\lambda_n}= +\infty$. Then the above sets $\mathscr{L}(I, M)$ are closed and nowhere dense in $\mathscr{S}^2(\lambda)$. 
\end{prop}
\begin{proof}
\proofpart{1}{Closed:} 
This part is simple, and only utilizes the fact that elements in $\mathscr{S}^2(\lambda)$ are integrable. Assume $(f_k)_k \subset \mathscr{L}(I, M)$ which converges to $f\in \mathscr{S}^2(\lambda)$. Fix an arbitrary $\eta>0$, and observe that for each $k$, we have
\begin{multline*}
m \left( \{ \zeta\in I: |f(\zeta)| > M+ \eta \} \right) \leq m \left( \{|f-f_k| > \eta  \} \right) + m \left( \zeta\in I: |f_k(\zeta)| \geq M \} \right) \\
\leq \frac{1}{\eta^2} \int_{\T} \abs{f-f_k}^2 dm \to 0, \qquad k\to \infty.
\end{multline*}
Letting $\eta \downarrow 0$ and using monotonicity of measures, we conclude that $m(I \cap \{|f|>M\})=0$, therefore it remains to show $m(I \cap \{|f|=M\})=0$. But this is even simpler and follows from the estimate
\begin{multline*}
m \left( I \cap \{ |f| = M \} \right) \leq \frac{2}{M^2} \int_{I} |f-f_k|^2 dm +  \frac{2}{M^2} \int_{I \cap \{ |f| = M \}} |f_k|^2 dm \\ = \frac{2}{M^2} \int_{I} |f-f_k|^2 dm + 2m(I \cap \{|f_k|=M\}) = \frac{2}{M^2} \int_{I} |f-f_k|^2 dm \to 0, \qquad k\to \infty.
\end{multline*}
We conclude that $f\in \mathscr{L}(I,M)$

\proofpart{2}{Nowhere dense:} Since $\mathscr{L}(I, M)$ is closed, it suffices to show that for any $g\in \mathscr{S}^2(\lambda)$ and any $\varepsilon_0>0$, there exists $g_{0} \in \mathscr{S}^2(\lambda) \setminus \mathscr{L}(I, M)$, such that 
\[
\norm{g-g_0}_{\mathscr{S}^2(\lambda)} < \varepsilon_0.
\]
By means of taking convolution with a smooth approximate of the identity, one can show that functions in $\mathscr{S}^2(\lambda)$ which are continuous on $\T$ are dense. Therefore, without loss of generality, we may assume that $g$ is continuous on $\T$. According to \thref{LEM:l2C}, we can exhibit an element $f_I \in \mathscr{S}^2(\lambda)$, with $\norm{f_I}_{\ell^2(\lambda)}=1$, which is essentially unbounded on $I$. It now readily follows that the elements 
\[
g_{\varepsilon,I}(\zeta)= g(\zeta) + \varepsilon \cdot f_I(\zeta), \qquad \varepsilon>0
\]
belong to $\mathscr{S}^2(\lambda)$, and are essentially unbounded on $I$ (here the continuity of $g$ ensures that it does not cancel the singularities of $f_I$). It follows that 
\[
\norm{g_{\varepsilon,I}- g}_{\mathscr{S}^2(\lambda)} = \varepsilon \norm{f_I}_{L^2} + \varepsilon < \varepsilon_0,
\]
provided that $\varepsilon>0$ is sufficiently small. We therefore conclude that $\mathscr{S}^2(\lambda) \setminus \mathscr{L}(I, M)$ is dense in $\mathscr{S}^2(\lambda)$. The proof is now complete.
\end{proof}

We demonstrate how our main result immediately follows from this observation.
\begin{proof}[Proof of \thref{THM:BAIRETOP}]
Pick a countable collection of arcs $(I_k)_k$ centered at a dense subset of points on $\T$, and with arbitrary small lengths. Pick also a sequence of positive numbers $(M_k)_k$ with $M_k \uparrow +\infty$, and note that according to \thref{PROP:BAIRE} the set
\[
\left\{ f\in \mathscr{S}^2(\lambda): f  \, \text{essentially bounded in a neighborhood of} \, \T \right\} = \bigcup_{k} \mathscr{L}(I_k, M_k),
\]
is of the first category in the Hilbert space $\mathscr{S}^2(\lambda)$. Applying the Baire category Theorem, we conclude that the subset of $\mathscr{S}^2(\lambda)$ consisting of functions which are unbounded in every neighborhood of $\T$ is generic in $\mathscr{S}^2(\lambda)$. This proves that $(i) \implies (iii) \implies (ii)$. But have already seen that the divergence of $\sum_n \frac{1}{\lambda_n} $ is necessary for $(ii)$ in view of the Cauchy-Schwarz inequality, hence the proof of our main result is complete.
\end{proof}

\section{Measures with small support}

\subsection{Regularity assumptions}
We begin this section by summarizing how the additional regularity assumption $(ii)$ of \thref{THM:BAIREPROB} will be implemented in this sequel. By means of substituting $w_n$ with $\min(w_n, 1)$, we may without loss of generality always assume that $w_n \leq 1$ for all $n$. The following lemma is an obvious modification of Lemma 2.1 in \cite{limani2025generic}, with $1/w_n$ playing the role of $\lambda_n$ therein. The simple details are thus left to the reader.

\begin{lemma}\thlabel{LEM:Regiiused} Let $(w_n)_{n=0}^\infty$ be positive numbers satisfying the hypothesis $(ii)$ in \thref{THM:BAIREPROB}. Then there exists a number $M(w)>1$ such that for any integer $M>M(w)$, there exists constants $A(M)>0$, such that the following statements hold:
\begin{enumerate}
\item[(a.)]
\[
\sum_{j=0}^n (1+j)^{M-1} w_j \leq A(M) n^M w_n, \qquad n=1,2,3,\dots
\]
\item[(b.)]
\[
\sum_{j\geq n} \frac{w_j}{j^{M+1}} \leq A(M) \frac{w_n}{n^M } , \qquad n=1,2,3,\dots
\]
\end{enumerate}
\end{lemma}

\subsection{Smooth functions with small amplitudes and small uniformly spread support}

The proof of \thref{THM:BAIREPROB} fundamentally relies on the construction of smooth functions with uniformly small amplitudes, and whose supports, though of arbitrarily small Lebesgue measure, are well-dispersed across $\T$. The following lemma may be viewed as an analogue to the so-called Ivashev-Musatov Lemma (for instance, see Lemma 3.2 in \cite{korner1986theorem3}, or P. 67 in \cite{havinbook}).

\begin{lemma}\thlabel{LEM:PRINCIPLEMMA}
For any integer $M>1$ there exists $c(M)>0$ such that the following statement holds: For any integer $N>1$, and for any $0<\delta<1$, we can find functions $\psi_{N,\delta} \in C^\infty(\T)$ satisfying the following properties:
\begin{enumerate}
    \item[(i)] $\psi_{N,\delta} \geq 0$ on $\T$,
    \item[(ii)] $\int_{\T} \psi_{N,\delta} dm =1$,
    \item[(iii)] for any arc $J\subset \T$ with $m(J)\geq \delta$, we have 
    \[
    m\left( \{\zeta \in J: \psi_{N,\delta}(\zeta)=0 \}\right) \geq 100^{-1} m(J),
    \]
     and the set $\{\zeta \in J: \psi_{N,\delta}(\zeta)=0 \}$ is a finite union of arcs,
    \item[(iv)]
    \[
    \abs{\widehat{\psi_{N,\delta}}(n)} \leq C(M)N^{-1/2} \min \left( 1, N^{-M} |n|^M, N^M |n|^{-M} \right), \qquad n\neq 0.
    \]
    \item[(v)] Every $\delta$-neighborhood of $\T$ meets $\supp{\psi_{N,\delta}}$.
\end{enumerate}
\end{lemma}

The proof of \thref{LEM:PRINCIPLEMMA} will be based on the following neat construction by T. W. K\"orner, which was used as a crucial building-block for proving Ivashev-Musatov type theorems in \cite{korner1986theorem3}, by-passing previous technical stationary phase arguments, also referred to as Van der Corput lemmas.

\begin{lemma}[T. W. K\"orner, Lemma 2.2 in \cite{korner1986theorem3}]\thlabel{LEM:KORNER22} For any integer $M>1$, there exists a constant $C(M)>0$, such that the following statement holds. For any closed arc $I\subseteq \T$ and for any integer $N\geq 1$, there exists functions $f_{I,N} \in C^\infty(\T)$ satisfying the following properties:
\begin{enumerate}
    \item[(a.)] $\supp{f_{I,N}}\subseteq I$,
    \item[(b.)] $-10\leq f_{I,N}\leq 10$ on $\T$,
    \item[(c.)] $\int_{I} \abs{f_{I,N}} dm = m(I)$,
    \item[(d.)] the set $\{\zeta \in I: f_{I,N}(\zeta)=- 10 \}$ is a finite union of arcs, which occupy a substantial proportion of $I$: $m\left( \{\zeta \in I: f_{I,N}(\zeta)= - 10 \} \right) \geq (50)^{-1}m(I)$.
    \item[(e.)] The following Fourier estimate holds:
    \[
    \abs{\widehat{f_{I,N}}(n)} \leq C(M) m(I) N^{-1/2} \min \left( 1, (|n|m(I)/N)^M, (|n|m(I)/N)^{-M} \right), \qquad n\neq 0.
    \]
\end{enumerate}
\end{lemma}
\begin{proof}[Proof of \thref{LEM:PRINCIPLEMMA}]
Fix arbitrary integers $M,N>1$, and let $0<\delta<1$ be arbitrary. 
Consider a finite partition $(I_k)_k$ of $\T$, consisting of closed arcs with mutually disjoint interiors, and a corresponding finite set of positive integers $(N_k)_k$, for which the following holds:
\[
\max_{k} m(I_k) \leq \delta/4, \qquad \min_k N_k \geq N, \qquad \frac{N_k}{m(I_k)} \, \text{non-increasing}.
\]
For instance, the constant sequence of integers $N_k = N(\delta) :=\max( 4/\delta, N)$ and $(I_k)_k$ to be a uniform partition of $\T$, consisting $N(\delta)$ closed arcs of length $N^{-1}(\delta)$, will do. We now apply \thref{LEM:KORNER22} in order to exhibit corresponding functions $f_k := f_{I_k, N_k}$, and form the sum 
\[
F(\zeta)= \sum_k f_{k}(\zeta), \qquad \zeta \in \T.
\]
Since $(I_k)_k$ have disjoint interior, we readily see from $(a.)$ of \thref{LEM:KORNER22} that $-10\leq F \leq 10$, while $(c.)$ implies that 
\[
\int_{\T} \abs{F} dm = \sum_k \int_{I_k} \abs{f_k} dm = \sum_k m(I_k) =1.
\]
We now claim that the desired functions are of the form:
\[
\psi_{N,\delta}(\zeta) := \frac{10+ F(\zeta)}{\int_{\T}(10+F)dm}, \qquad \zeta \in \T.
\]
Indeed, the properties $(i)$ and $(ii)$ readily follow by construction. To verify $(iii)$, we fix an arbitrary arc $J\subset \T$ with $m(J)\geq \delta$. Since $m(I_k)\leq \delta/4$ for all $k$, we can find a finite sub-collection $(I'_j)_j$ so that the union $\cup_k I'_j$ is contained in $J$, and $\sum_k m(I'_k) \geq m(J)/2$. Recall that the property $(c.)$ of \thref{LEM:KORNER22} implies that $\{\zeta \in \T: \psi_{N,\delta}(\zeta)\} = \{ F=-10 \}$, is again a finite union of closed arcs, hence we get 
\begin{multline*}
m(\{\zeta \in \T: \psi_{N,\delta}(\zeta)\}) \geq m(\{ \zeta \in \cup_{j}I'_j : F(\zeta)=-10 \}) \\
= \sum_j m(\{\zeta \in I'_j: f_j(\zeta)=-10 \}) \geq \frac{1}{50}\sum_j m(I'_j) \geq \frac{1}{100} m(J).
\end{multline*}
This proves the claim in $(iii)$. In order to verify $(v)$, we note that $\{\zeta \in \T: \psi_{N,\delta} >0 \}= \cup_k \{\zeta \in I_k: f_k(\zeta)>-10 \}$, where each set $\{\zeta \in I_k: f_k(\zeta)>-10 \}$ have non-trivial interior. Now since $(I_k)_k$ forms a partition of $\T$ and $\max_k m(I_k) \leq \delta/4$, it readily follows that $\sup_{\zeta\in \T}\dist{\zeta}{\supp{\psi_{N,\delta}}}\leq \delta$, therefore $(v)$ holds. It only remains to prove $(iv)$. To this end, we primarily note that since $F$ is real-valued and the estimate
\[
\int_{\T}(10+F)dm \geq 10 - \int_{\T}\abs{F} dm \geq 9,
\]
holds, we see that $\abs{\widehat{\psi_{N,\delta}}(n)} \leq \abs{\widehat{F}(n)}$ for $n\neq 0$, hence it suffices to prove the estimate for the positive Fourier coefficients of $F$. Using property $(e.)$ for each $f_k$, and the fact that they have disjoint support, we easily see that 
\begin{multline*}
\abs{\widehat{F}(n)} \leq \sum_k \abs{\widehat{f_k}(n)} \leq C(M) \sum_k m(I_k) N^{-1/2}_k \min (1, (n m(I_k)/N_k)^M) \\
\leq  C(M) N^{-1/2} \min (1, (n/N)^M) \sum_k m(I_k) = C(M) N^{-1/2} \min (1, (n/N)^M),
\end{multline*}
where we also used the fact that $\inf_k N_K \geq N$. Now, utilizing the property that $N_k/m(I_k)$ was chosen to be non-increasing (possibly constant), we also have 
\[
\abs{\widehat{F}(n)} \leq \sum_k \abs{\widehat{f_k}(n)} \leq C(M) N^{-1/2} \sum_k m(I_k) (nm(I_k)/N_k)^{-M} \leq C(M) N^{-1/2} (n/N)^{-M}.
\]
This is enough to conclude that $(iv)$ holds, hence the proof of the lemma is complete.
\end{proof}

\subsection{Another Baire category argument} 
Before carrying out the main Baire category argument, we shall need to properly set up the appropriate functional analytic setting. The following basic lemma, whose proof we only sketch, contains these necessary ingredients.

\begin{lemma}\thlabel{LEM:FUNCAN} Let $w=(w_n)_n$ be any sequence of positive numbers, and denote by $\mathscr{P}_{w}$ the collection of order pairs $(\sigma,E) \in \ell^2(w) \times \mathscr{C}$ satisfying the properties
\begin{enumerate}
    \item[(a.)] $\supp{\sigma} \subseteq E$,
    \item[(b.)] $\sigma$ is a positive finite Borel measure on $\T$.
\end{enumerate}
Then the set $\mathscr{P}_{w}$ equipped with the metric
\[
d_{w}\left( (\sigma,E), (\tau,F) \right) = d_{\mathscr{C}}(E,F) + \norm{\sigma-\tau}_{\ell^2(w)},
\] 
forms a complete metric space $(\mathscr{P}_w, d_w)$.
\end{lemma}
\begin{proof}
Note that since the metric spaces $(\mathscr{C}, d_{\mathscr{C}})$ and $(\ell^2(w), \norm{\cdot}_{\ell^2(w)})$ are both complete, it readily follows that $\ell^2(w) \times \mathscr{C}$ is also complete in the metric $d_w$, hence it remains only to verify that the sub-collection $\mathscr{P}_w$ is sequentially closed in $\ell^2(w) \times \mathscr{C}$. Suppose $\{(\sigma_n, E_n)\}_n \in \mathscr{P}_w$ with  $d_w(\sigma_n, E_n),(\sigma, E)) \to 0$, where $(\sigma,E) \in \ell^2(w) \times \mathscr{C}$. Now for any $f\in C^\infty(\T)$ with $f\geq 0$ on $\T$, we have 
\[
\sum_m \widehat{\sigma}(m) \widehat{f}(m) = \lim_n \sum_m \widehat{\sigma_n}(m) \widehat{f}(m)= \lim_n \int_{\T} f(\zeta) d\sigma_n(\zeta) \geq 0.
\]
It follows from this calculation that $\sigma$ is a positive finite Borel measure on $\T$, hence we only need to show that $\supp{\sigma} \subseteq E$. Now if $U$ is an arbitrary neighborhood of $\T \setminus E$, then since $d_{\mathscr{C}}(E_n,E) \to 0$, there exists $N=N(U)>0$, such that $U \cap E_n = \emptyset$ for all $n> N$. Using the assumption that $\supp{\sigma_n} \subseteq E_n$ for each n, we see that 
\[
\int_{\T} f(\zeta) d\sigma(\zeta) = \sum_m \widehat{\sigma}(m) \widehat{f}(m) = \lim_n \sum_m \widehat{\sigma_n}(m) \widehat{f}(m)= \lim_n \int_{\T} f(\zeta) d\sigma_n(\zeta) =0,
\]
for any $f\in C^\infty(\T)$ that vanishes identically off the set $U$. This proves the claim.

\end{proof}

Stacking up the results from the previous subsection, we now come to the crucial step in the Baire category argument. Here, the regularity assumptions on $w$ will be utilized in essential ways.

\begin{prop}\thlabel{PROP:Ea} Let $(w_n)_n$ be positive real numbers, satisfying the hypothesis $(i)-(ii)$ of \thref{THM:BAIREPROB}. For any number $0<\gamma<1$, the collection of ordered pairs
\[
\mathscr{E}_{\gamma} := \left\{(\sigma,E)\in \mathscr{P}_w: \, \, m(E)< \gamma \right\},
\]
is open and dense in the complete metric space $(\mathscr{P}_w, \rho_w)$.

\end{prop}

The proof is a bit long, and we have for the readers convenience divided it into different steps. 

\begin{proof}
    We have already seen that $(\mathscr{P}_w, \rho_w)$ forms a complete metric space in \thref{LEM:FUNCAN}, thus we proceed to verify the statements on openness and density. Below, we shall carry out the proof under the second assumption in $(ii)$, that there exists $C>1$ such that $C^{-1}w_n \leq w_k \leq C w_n$, for all $n\leq k \leq 2n$. The simple adaptions for the case $w_n \downarrow 0$ will be outlined right after.
    
    \proofpart{1}{Openness:} Let $0<\gamma<1$, fix an arbitrary element $(\sigma,E) \in \mathscr{E}_\gamma$, which at its turn allows us pick a small number $0<\delta < \gamma - m(E)$. By the definition of Lebesgue measure, we can cover $E$ by a sequence of open arcs $(I_j)_j$ with each $m(I_j)\leq \delta$, and $\sum_j m(I_j) \leq m(E)+ \delta$. Invoking compactness of $E$, we can pass to a finite sub-collection, say $I_1, I_2, \dots, I_l$, again with the same properties. Since we only have a finite sub-collection, we can define the positive number $\delta_0 := \min_{1\leq j\leq l} m(I_j)>0$. Now, note that the arcs $I_1, I_2, \dots, I_l$ also cover $K$, whenever $d_{\mathscr{C}}(E,K) < \delta_0/2$, hence it follows that $m(K) \leq \sum_j m(I_j) \leq m(E) + \delta < \gamma$. This implies that $\mathscr{E}_\gamma$ is open.
    
    \proofpart{2}{Reduction to smooth functions:}
    To prove that $\mathscr{E}_\gamma$ is dense, it suffices to show that for any $(\sigma,E) \in \mathscr{P}_w$ and any $\varepsilon>0$, there exists $(f,K) \in \mathscr{E}_\gamma$ with $\rho_w\left( (\sigma,E), (f, K) \right) < \varepsilon$. Now, for any $\varepsilon>0$, we can form $C^\infty$-smooth kernels $k_{\varepsilon}$ with the properties:
\[
k_{\varepsilon}\geq 0 \,\, \, \text{on} \, \, \, \T, \qquad \int_{\T} k_{\varepsilon} dm=1, \qquad \supp{k_{\varepsilon}}\subseteq I_{\varepsilon},
\]
where $I_{\varepsilon}$ is a $\varepsilon$-neighborhood of $\zeta=1$. If we set $f_{\varepsilon}:= \sigma \ast k_{\varepsilon}$, it follows that $\supp{f_{\varepsilon}} \subseteq K + I_{\varepsilon}:=\{\zeta \in \T: \dist{\zeta}{K}\leq \varepsilon \}$, whence $d_{\mathscr{C}}(\supp{f_{\varepsilon}}, E) \leq \varepsilon$, and therefore the ordered pair $(f_{\varepsilon}, \supp{f_{\varepsilon}})$ belongs to $\mathscr{P}_w$. Furthermore, since
\[
\abs{\widehat{f_{\varepsilon}}(n)-\widehat{\sigma}(n)} =\abs{\widehat{\sigma}(n)}\cdot \abs{\widehat{k_\varepsilon}(n)-1} \leq 2\abs{\widehat{\sigma}(n)} , \qquad n\in \mathbb{Z},
\]
it follows that for any $\varepsilon>0$, and any integer $N>0$, that
\[
\sum_n \abs{\widehat{f_{\varepsilon}}(n)-\widehat{\sigma}(n)}^2 w_{|n|} \lesssim \sum_{\abs{n}\leq N}\abs{\widehat{\sigma}(n)}^2\cdot \abs{\widehat{k_\varepsilon}(n)-1}^2 w_{|n|} + \sum_{|n|>N} \abs{\widehat{\sigma}(n)}^2 w_{|n|}.
\]
Sending $\varepsilon \to 0+$ first, then $N \to \infty$, we conclude that $f_{\varepsilon} \to \sigma$ in $\ell^2(w)$. We draw the conclusion that by means of passing to such a pair $(f_{\varepsilon}, \supp{f_{\varepsilon}})$, we may without loss of generality assume that $\sigma \in C^\infty(\T)$. 

\proofpart{3}{Estimations using the principal lemma:}

Moving forward, we now make use of the crucial assumption $(i)$ in \thref{THM:BAIREPROB}, which allows us to exhibit positive integers $(N_j)_j$ for which $w_{N_j} \to 0$. Appealing to \thref{LEM:PRINCIPLEMMA} with parameters $M=M(w)>0$ as in \thref{LEM:Regiiused}, and $(N_j)_j$ as in the previous paragraph, we can exhibit corresponding functions $\psi_{j}:= \psi_{N_j, w_{N_j} } \in C^\infty(\T)$ which satisfy the following properties:
\begin{enumerate}
    \item[(i)] $\psi_{j} \geq 0$ on $\T$,
    \item[(ii)] $\int_{\T} \psi_{j} dm =1$,
    \item[(iii)] for any arc $J\subset \T$ with $m(J)\geq w_{N_j}$, we have 
    \[
    m\left( \{\zeta \in J: \psi_{j}(\zeta)=0 \}\right) \geq 100^{-1} m(J),
    \]
    and the set $\{\zeta \in J: \psi_{j}(\zeta)=0 \}$ is a finite union of arcs,
    \item[(iv)]
    \[
    \abs{\widehat{\psi_{j}}(n)} \leq C(M)N_j^{-1/2} \min \left( 1, N_j^{-M} |n|^M, N_j^M |n|^{-M} \right), \qquad n\neq 0.
    \]
    \item[(v)] Every $w_{N_j}$-neighborhood of $\T$ meets $\supp{\psi_j}$.
 \end{enumerate}
  Now set $f_j := \sigma \cdot \psi_{j}$, and let $E_j := E \setminus U_{j}$, where $U_{j}$ are the interiors of the arcs comprising the set $\left\{\zeta \in \T: \psi_j(\zeta)=0 \right\}$. It follows from $(i)$ and $(iii)$ $(f_j,E_j) \in \mathscr{P}_w$, while the property $(v)$ ensures that $d_{\mathscr{C}}(E,E_j)\leq w_{N_j}$. We now show that the property $(iii)$ implies that 
  \begin{equation}\label{EQ:kappa}
  m(E_j) \leq (1-10^{-4}) m(E),
   \end{equation}
  whenever $j$ sufficiently large. Since $E\subset \T$ is compact with $m(E)>0$ (it supports a non-trivial smooth function $\sigma$), we can find a finite collection of closed arcs $I_1, \dots, I_l$, such that 
\[
E \subseteq \cup_{k=1}^l I_k, \qquad \frac{101}{100} m(E) \geq \sum_{k=1}^l m(I_k).
\]
But then by definition of $E_j$, if me choose $j$ such that $0<w_{N_j} < \min_{1\leq k\leq l} m(I_k)$, then it follows from property $(iii)$ that
\[
m(E_j ) \leq \sum_{k} m(\{\zeta \in I_k : \psi_j(\zeta)>0 \}) \leq (1-\frac{1}{100}) \sum_{k} m(I_k) \leq (1-\frac{1}{100}) (1+\frac{1}{100}) m(E).
\]
This proves \eqref{EQ:kappa}. Next, we proceed by showing that 
\[
\sum_{n} \abs{\widehat{f_j}(n)-\widehat{\sigma}(n)}^2w_{|n|} \to 0, \qquad j \to \infty.
\]
Since all functions involved are real-valued, we only need to estimate the sum for $n\geq 0$. Since the parameter $M=M(w)>0$ and the function $\sigma$ have both been fixed, we shall for the sake of clarity, disregard constants depending on $M$ and $\sigma$, and simply write $\lesssim$ to indicate inequalities up to such constants. First, we note that the property $(ii)$ of $\psi_j$ allows us to write:
\[
\widehat{f_j}(n)-\widehat{\sigma}(n) = \sum_{m\neq n} \widehat{\sigma}(m) \widehat{\psi_j}(n-m) = \left(\sum_{|m|\leq n/2} + \sum_{\substack{|m|>n/2 \\ m\neq n}} \right) \widehat{\sigma}(m) \widehat{\psi_j}(n-m).
\]
Using the property $(iii)$ and the assumption that $\sigma$ is smooth on $\T$, we get 
\[
\abs{\sum_{\substack{|m|>n/2 \\ m\neq n}} \widehat{\sigma}(m) \widehat{\psi_j}(n-m)} \lesssim N^{-1/2}_j \sum_{|m|>n/2} \abs{\widehat{\sigma}(m)} \lesssim   \frac{1}{N^{1/2}_j(1+n)}.
\]
We obtain
\begin{equation}\label{EQ:Largem}
\sum_{n\geq 0} \abs{\sum_{\substack{|m|>n/2 \\ m\neq n}} \widehat{\sigma}(m) \widehat{\psi_j}(n-m)}^2 w_{n} \lesssim \frac{1}{N_j} \sum_{n\geq 0} \frac{1}{(1+n)^2} \lesssim \frac{1}{N_j}.
\end{equation}
Thus it remains to estimate the sum
\[
\sum_{n\geq 0} \abs{ \sum_{|m|\leq n/2} \widehat{\sigma}(m) \widehat{\psi_j}(n-m)}^2 w_{n}.
\]
To this end, we shall split the above sum in the following way:

\begin{equation}\label{EQ:Nepssum}
\sum_{n\geq 0} \abs{ \sum_{|m|\leq n/2} \widehat{\sigma}(m) \widehat{\psi_j}(n-m)}^2 w_{n} = \left(\sum_{0\leq n < N_j} + \sum_{n\geq N_j} \right) \abs{ \sum_{|m|\leq n/2} \widehat{\sigma}(m) \widehat{\psi_j}(n-m)}^2 w_{n}.
\end{equation}
For sufficiently large $n\geq N_j$, we again use property $(iii)$ of $\psi_{j}$ in conjunction with $(b.)$ of \thref{LEM:Regiiused}, which yields
\begin{multline*}
\sum_{n\geq N_j} \abs{ \sum_{|m|\leq n/2} \widehat{\sigma}(m) \widehat{\psi_j}(n-m)}^2 w_{n} \leq \sum_{n\geq N_j} \left( \sum_{|m|\leq n/2} \abs{\widehat{\sigma}(m) \widehat{\psi_j}(n-m)} \right)^2 w_{n}  \\ 
\lesssim N_j^{2M-1} \sum_{n\geq N_j}\left( \sum_{|m|\leq n/2} \frac{\abs{\widehat{\sigma}(m)}}{|n-m|^{M}} \right)^2 w_n \\
\lesssim N_j^{2M-1} \sum_{n\geq N_j} \frac{w_n}{n^{2M}} \left( \sum_{|m|\leq n/2} \abs{\widehat{\sigma}(m)} \right)^2 \lesssim \norm{\widehat{\sigma}}^2_{\ell^1} w_{N_j}.
\end{multline*}
On the other hand, when $0\leq n<N_j$ we use the assumption $(iii)$ of $\psi_j$ with $(a.)$ of \thref{LEM:Regiiused}, in order to deduce
\begin{multline}\label{EQ:nestiii}
\sum_{0\leq n < N_j} \abs{\sum_{|m|\leq n/2} \widehat{\sigma}(m) \widehat{\psi_j}(n-m)}^2 w_{n} \leq 
\sum_{0\leq n < N_j} \left(\sum_{|m|\leq n/2} \abs{\widehat{\sigma}(m)} \abs{ \widehat{\psi_j}(n-m)} \right)^2 w_{n} \\
\lesssim \frac{1}{N^{2M+1}_j} \sum_{0<\leq n< N_j} w_n \left( \sum_{|m|<n/2} \abs{\widehat{\sigma}(m)}|n-m|^M \right)^2 \lesssim \frac{1}{N^{2M+1}_j} \sum_{0\leq n< N_j} (1+n)^{2M} w_n \norm{\widehat{\sigma}}^2_{\ell^1} \\
\lesssim  \norm{\widehat{\sigma}}^2_{\ell^1} w_{N_j}.
\end{multline}
Combining the estimates in \eqref{EQ:Largem} and in \eqref{EQ:Nepssum}, we arrive at 
\[
\sum_n \abs{\widehat{f_j}(n)-\widehat{\sigma}(n)}^2 w_{|n|} \leq C(M,\sigma) \max (N^{-1}_j , w_{N_j}),
\]
where $C(M,\sigma)>0$ is a constant only depending on $M$ and $\sigma$. Recalling that $w_{N_j} \to 0$ and \eqref{EQ:kappa}, we have, upon on a standard re-scaling argument, established the following statement. For any pair $(\sigma,E) \in \mathscr{P}_w$ with $\sigma \in C^\infty(\T)$, and for any $\varepsilon>0$ sufficiently small, there exists functions $f_{\varepsilon} \in C^\infty(\T)$, which satisfy 
\begin{equation}\label{EQ:e_j}
\norm{f_{\varepsilon}\cdot \sigma-\sigma}_{\ell^2(w)} \leq \varepsilon, \qquad m(E_{\varepsilon})\leq (1- 10^{-4}) m(E), \qquad \sup_{\zeta\in \T} \dist{\zeta}{E_{\varepsilon}} \leq \varepsilon,
\end{equation}
where $E_{\varepsilon} := E \setminus U_{\varepsilon}$, and $U_{\varepsilon}$ denotes the finite union of open arcs comprising the interior of the set $\{\zeta \in \T: f_{\varepsilon}(\zeta)=0 \}$.

\proofpart{4}{An iterative scheme:} According to \eqref{EQ:e_j}, we can produce elements $(f_\varepsilon, E_\varepsilon) \in \mathscr{P}_w$ which are close to $(\sigma,E)$ in the $\rho_w$-metric. Indeed, $E_\varepsilon \subseteq E$, we have $d_{\mathscr{C}}(E,E_\varepsilon)\leq \sup_{\zeta\in \T} \dist{\zeta}{E_{\varepsilon}} \leq \varepsilon$. However, we have only been able to contract the Lebesgue measure of $E_\varepsilon$ by a factor $(1-10^{-4})$ to that of $E$. In order to ensure that our approximates $(f_\varepsilon, E_\varepsilon)$ of $(\sigma,E)$ belong to $\mathscr{E}_\gamma$, we shall recursively implement the procedure in \eqref{EQ:e_j}. To this end, pick an integer $L=L(\gamma)>1$, such that $(1- 10^{-4})^L \leq \gamma$. Fix an arbitrary $\varepsilon >0$, and set $(f_0, E_0)= (\sigma,E)$ and applying the previous argument with $\varepsilon_1 := \varepsilon/L$, we can exhibit a pair $(f_{\varepsilon_1},E_1)$ satisfying the properties in \eqref{EQ:e_j}. We now set $f_1 := f_{\varepsilon_1} \cdot \sigma$ and apply \eqref{EQ:e_j} with $(f_1, E_1)$ in place of $(f_0, E_0)$. Continuing in this manner, we obtain at the $L$-th stage, a smooth function $f_{\varepsilon_L}$ on $\T$ and a corresponding set $U_{\varepsilon_L}$, being the interior of arcs comprising the set $\{\zeta\in \T: f_{\varepsilon_L}(\zeta)=0 \}$, such that 
\[
\rho_w \left( (f_{\varepsilon_L}\cdot f_{L-1}, E_{L}), (f_{L-1},E_{L-1}) \right) \leq \varepsilon_L \qquad m(E_{\varepsilon_L}) \leq (1- 10^{-4}) m(E_{L-1}),
\]
with $E_L := E_{L-1} \setminus U_{\varepsilon_L}$. Setting $f_L := f_{\varepsilon_L}f_{L-1} = \sigma \cdot \prod_{j=1}^L f_{\varepsilon_j}$, we obtain 
\[
\rho_w \left( (f_L, E_L), (\sigma,E) \right) \leq \sum_{j=1}^L \rho_w \left( (f_{j}, E_{j}), (f_{j-1},E_{j-1}) \right) \leq \varepsilon,
\]
and $m(E_L) \leq (1- 10^{-4}) m(E_{L-1}) \leq (1- 10^{-4})^L m(E) \leq \gamma$, thus the pair $(f_L, E_L) \in \mathscr{E}_\gamma$. We therefore conclude that $\mathscr{E}_\gamma$ is dense in $\mathscr{P}_w$, for any fixed $\gamma>0$. The proof is now complete.
\end{proof}

\begin{proof}[Proof of \thref{THM:BAIREPROB}]
Now pick a sequence $(\gamma_n)_n$ with $\gamma_n \downarrow 0$ and consider the set $\mathscr{E} := \cap_n \mathscr{E}_{\gamma_n}$, where each $\mathscr{E}_{\gamma_n}$ as in \thref{PROP:Ea}. According to the Baire category Theorem, we conclude that $\mathscr{E}$ is dense in the complete metric space $(\mathscr{P}_w, \rho_w)$. Now whenever $(\sigma,E) \in \mathscr{E}$, then $\supp{\sigma}\subseteq E$ with $m(E)\leq \gamma_n$, for all $n$. It readily follows that $\supp{\sigma}$ has Lebesgue measure zero in $\T$.
\end{proof}

We complete this section by demonstrating how \thref{THM:BAIREPROB} is proved when $w_n \downarrow 0$. As the proof is very similar, we only indicate the distinct features. First, we note that property $(b.)$ in \thref{LEM:Regiiused} trivially holds for monotonic $(w_n)_n$, while $(a.)$ may fail. However, the property $(a.)$ was only utilized in combination with the estimate
\[
\abs{\widehat{\psi_{N}}(n)}\leq c(M) N^{-M-1/2} |n|^M, \qquad n\neq 0,
\]
as seen from \eqref{EQ:nestiii}. Now if we instead utilize the simpler estimate $\abs{\widehat{\psi_{N}}(n)}\leq c(M) N^{-1/2}$ in \eqref{EQ:nestiii}, we obtain 
\[
\sum_{0\leq n < N} \abs{\sum_{|m|\leq n/2} \widehat{\sigma}(m) \widehat{\psi_N}(n-m)}^2 w_{n} \lesssim \norm{\widehat{\sigma}}^2_{\ell^1} N^{-1} \sum_{0\leq n< N} w_n,
\]
where the left hand side also tends to zero, since it is a Cesar\'o mean of $(w_n)_n$. The rest of the proof remains unchanged.

\subsection{Measures with eventually sparse Fourier spectrum}
This subsection is devoted to proving that the assumption $\liminf_n  w_n =0$ alone is not enough to ensure the existence of elements with arbitrary small support in the corresponding space $\ell^2(w)$. The result will principally be based on \thref{LEM:UPZYG}, which is essentially an uncertainty principle based on the work of A. Zygmund on lacunary series. Below, shall further clarify this point. Given a number $N>1$, we declare a subset of integers $\Lambda \subset \mathbb{Z}$ to be $N$-sparse, if 
\[
\inf_{\substack{n,m\in \Lambda \\ |n|>|m|}} \frac{|n|}{|m|}\geq N.
\] 
By default, we use the convention that the empty set is $N$-sparse for any $N>1$. The following result is likely well-known to experts, and asserts that measures which eventually have their Fourier spectrum contained in a $3$-sparse set must have full support in $\T$. Since we have not been able to locate a proof of this in the literature, we shall for the readers convenience include a brief sketch of proof.


\begin{lemma}[A. Zygmund] \thlabel{LEM:ZYGlac} Let $\Lambda \subset \mathbb{Z}$ be a $3$-sparse set. Then whenever $\mu$ is complex finite Borel measure on $\T$ with following properties:
\begin{enumerate}
    \item[(a.)] There exists $N_0>0$, such that $\supp{\widehat{\mu}} \cap \{|n|>N_0\} \subseteq \Lambda$
    \item[(b.)] there exists an arc $I \subset \T$ such that $\mu(I)=0$,
\end{enumerate}
then $\mu\equiv 0$.
\end{lemma}
\begin{proof}
Note that by a simple matter of translation, we may assume that $\mu$ vanishes on the arc $I_\eta$ centered at $\zeta=1$ and of length $\eta>0$. Now set $\Lambda_{0} := \Lambda \cap \{|n|>N_0\}$ and note that $\Lambda_0$ is also $3$-sparse, hence there exists a numerical constant $c_0>0$, such that for any $\lambda_j \in \Lambda_0$, the sets
\[
S_j := \{ n \in \mathbb{Z}: 1\leq |n-\lambda_j| \leq c_0 \abs{\lambda_j} \},
\]
do not intersect $\Lambda_0$. Now, since the Fourier spectrum of $\mu$ is eventually contained in $\Lambda_0$ (that is, whenever $|n|>N_0$), we have
\[
\widehat{\mu}(\lambda_j) = \sum_k \widehat{\mu}(\lambda_j -k ) \widehat{F_{N_j}}(k) =\int_{\T \setminus I_{\eta}} \zeta^{-\lambda_j} F_{N_j}(\zeta) d\mu(\zeta), 
\]
where $F_N$ denotes the classical Fej\'er kernel, and say, $c_0 \lambda_j/10 \leq N_j < c_0 \lambda_j$ is an integer. Now, using the trivial point-wise bound $\abs{F_N(e^{it})}\leq \frac{1}{N \sin^2(t/2)}$, we get 
\[
\abs{\widehat{\mu}(\lambda_j)} \leq 100 \eta^{-2} \frac{1}{N_j} \leq \frac{1000}{c_0 \eta^2} \lambda^{-1}_j,
\]
for all $\lambda_j \in \Lambda_0$ with $\abs{\lambda_j}\geq N_1$, say. Now since $\Lambda_0$ is $3$-sparse, we conclude that there $(\widehat{\mu}(n))_n \in \ell^1$, hence there exists a continuous function $f$ on $\T$ with $d\mu = f dm$. Now using an argument similar to that of \thref{LEM:RIGIDLAC} (see also exercise 1.4 in \cite{katznelson2004introduction}), one can show that there exists a constant $C(I_\eta, \Lambda)>0$, such that 
\[
\sum_{|j|>j_0} \abs{\widehat{f}(\lambda_j)} \leq C(I_\eta, \Lambda) \sup_{\zeta \in I_\eta} \abs{f(\zeta)} =0.
\]
for large enough $j_0=j_0(\eta, \Lambda)>0$. We conclude that $f$ is a trigonometric polynomial which vanishes identically on $I_\eta$, hence $f\equiv 0$, and thus $\mu \equiv 0$.
\end{proof}

With this result at hand, we now establish the uncertainty principle which lays our foundation for \thref{THM:NECREG}.
\begin{proof}[Proof of \thref{LEM:UPZYG}] The proof is inspired from the K\"orner-Meyer Lemma (see Lemma 2.1 in \cite{korner2003topological}), and is also proved by contradiction. Suppose the statement fails for a certain set of parameters $(N, \gamma, \delta)$. Then there exists numbers $(\varepsilon_k)_k$ with $\varepsilon_k \to 0$, such that the following statements hold: there is a sequence of complex finite Borel measures $(\mu_k)_k$ on $\T$ with $\norm{\mu_k}\leq 1$, and 
\[
(i)_k \, \sum_{|n|\leq N} \abs{\widehat{\mu_k}(n)}^2 \geq \gamma, \qquad (ii)_k \, \abs{\widehat{\mu_k}(n)} \leq \varepsilon_k \, \, \text{for} \, |n|>N \, \text{off a 3-sparse set}
\]
for all $k$, but yet $\sup_{\zeta\in \T}\dist{\zeta}{\supp{\mu_k}}> \delta$. Now, consider a uniform partition $\{J_1, J_2, \dots, J_l\}$ of $\T$ consisting of arcs of equal lengths, not exceeding $\delta/10$. By the Pigeon hole principle, there must exist an arc $J_{l_0}$ among these, for which the support of infinitely many of the measures $(\mu_k)_k$ do not meet. This implies that we can find an infinite sequence $(k_j)_j$ for which $\norm{\mu_{k_j}} \leq 1$ and the following statements hold:
\[
\supp{\mu_{k_j}} \cap J_{l_0} = \emptyset, \qquad \sum_{|n|\leq N} \abs{\widehat{\mu_{k_j}}(n)}^2 \geq \gamma, \qquad \, \abs{\widehat{\mu_{k_j}}(n)} \leq \varepsilon_{k_j} \, \, \text{for} \, |n|>N \, \text{off a 3-sparse set}.
\]
Invoking Helly's selection Theorem, we may extract a weak-star convergent cluster point $\mu$, which again is a complex finite Borel measure on $\T$, that satisfies the following properties:
\[
\norm{\mu}\leq 1, \qquad \supp{\mu} \cap J_{l_0} = \emptyset, \qquad \sum_{|n|\leq N} \abs{\widehat{\mu}(n)}^2 \geq \gamma, \qquad \supp{\widehat{\mu}} \cap \{|n|>N\} \, \, \text{is 3-sparse}.
\]
According to \thref{LEM:ZYGlac}, the second and the forth properties of $\mu$ forces $\mu \equiv 0$, which is in obvious violation to the third property. We have thus reached a desired reductio ad absurdum, and the claim follows.
\end{proof}

Finally, we may now complete our main task in this subsection.

\begin{proof}[Proof of \thref{THM:NECREG}] We shall construct the desired sequence $(w_n)_n$ by means of running \thref{LEM:UPZYG} inductively. To this end, set $\varepsilon_1=N_1=1$ and suppose by means of induction that $\varepsilon_1 > \dots > \varepsilon_k>0$ and $N_1 < \dots < N_k$ have already been constructed. Pick a number $N_{k+1} > 3N_k$ and apply \thref{LEM:UPZYG}, which implies that there exists a number $0<\varepsilon_{k+1} < \varepsilon_k/2$, such that whenever $\mu$ is a complex finite Borel measure on $\T$ with $\norm{\mu}\leq 1$ and satisfies the properties:
\begin{equation*}
(i)_{k} \, \sum_{|n|\leq N_{k+1}} \abs{\widehat{\mu}(n)}^2 \geq \varepsilon_k, \qquad (ii)_k \, \abs{\widehat{\mu}(n)} \leq \varepsilon_{k+1} \, \, \text{for} \, \, |n|>N_{k+1} \, \, \text{off a 3-sparse set},
\end{equation*}
then $\sup_{\zeta\in \T} \dist{\zeta}{\supp{\mu}}\leq 1/k$. This gives rise to a pair of sequences $(\varepsilon_k)_k, (N_k)_k$ with $\varepsilon_k \downarrow 0$ and $N_k \uparrow \infty$, which we now utilize in order to define $(w_n)_n$ as follows:
\[
w_n = 
\begin{cases}
    \varepsilon_k & \, \text{if} \, \, n=N_k, \\
    \frac{1}{\varepsilon^2_k} & \, \text{if} \, \, N_k <n <N_{k+1},
\end{cases} \qquad k=1,2,3,\dots
\]
Clearly $\liminf_n w_n =0$. Now let $\mu$ be a non-trivial complex finite Borel measure on $\T$ with $(\widehat{\mu}(n))_n \in \ell^2(w)$, which we may assume also satisfies $\norm{\mu} \leq 1$. We must verify that $(i)_k$ and $(ii)_k$ hold for all sufficiently large $k$. Now $(i)_k$ readily follows from the fact the partial $\ell^2$-sums of $\mu$ must eventually become strictly positive (possibly even infinite), while the right hand side tends to zero. Hence there exists $k_0>0$ such that $(i)_k$ holds for all $k>k_0$. In order to verify $(ii)_k$, we recall that since $\inf_k N_{k+1}/N_k \geq 3$ by construction, hence it suffices to show that 
\[
\abs{\widehat{\mu}(n)} \leq \varepsilon_{k+1}, \qquad N_{k+1}< |n| < N_{k+2},
\]
for large enough $k$. Since the Fourier coefficients of $\mu$ belong to $\ell^2(w)$, we can by definition of $(w_n)_n$ find an integer $k_1>0$, such that 
\[
\sum_{k>k_1} \frac{1}{\varepsilon^2_k} \sum_{N_k < |n|<N_{k+1}} \abs{\widehat{\mu}(n)}^2 = \sum_{k>k_1} \sum_{N_k < |n|< N_{k+1}} \abs{\widehat{\mu}(n)}^2 w_{|n|} \leq 1.
\]
But then it readily follows that
\[
\sup_{N_{k}<|n|<N_{k+1}} \abs{\widehat{\mu}(n)}^2 \leq \varepsilon^2_k \sum_{N_k < |n| < N_{k+1}} \abs{\widehat{\mu}(n)}^2 w_{|n|} \leq \varepsilon^2_k, \qquad k>k_1.
\]
Since $(\varepsilon_k)_k$ is monotonically decreasing, we have 
\[
\sup_{|n| \in \cup_{j\geq k} (N_j, N_{j+1})} \abs{\widehat{\mu}(n)} \leq \varepsilon_k,
\]
hence $(ii)_k$ holds for all $k>k_1$ large enough. We have thus verified that $(i)_k-(ii)_k$ hold for all $k>k_0+k_1$, and therefore $\supp{\mu}= \T$.
\end{proof}
As indicated in Section 2, the proof shows that one has complete freedom in modifying the sequence $(w_n)_n$ along any $3$-sparse subsequence, as long as it tends to zero, and it may decay there at any prescribed rate. However, this flexibility is compensated by the requirement that $(w_n)_n$ must be sufficiently large outside that subsequence, and the necessary size is implicitly determined by the uncertainty principle in \thref{LEM:UPZYG}, which is qualitative in nature, and thus does not yield explicit quantitative bounds.

\subsection{A simple construction involving classical Riesz products}
At last, we shall use the notion of the classical family of Riesz products introduced in the previous section, in order to construct a probability measure supported on a set of Lebesgue measure zero with Fourier coefficients in $\ell^2(w)$, under a slightly milder regularity assumption on $w$. We shall prove the existence under the assumption that $w_n$ tends to zero along an $\varepsilon N_j$-neighborhood of integers, centered at a very lacunary sequence $N_j$, which is arguably much stronger than merely requiring $\liminf_n w_n =0$.

\begin{thm}\thlabel{THM:singl2w} Let $(w_n)_n$ positive real numbers with the property that there exists a number $\varepsilon>0$ and positive integers $(N_j)_j$, such that
\[
W_j:= \sup_{n: \, |n -N_j|\leq \varepsilon N_j }  w_n \to 0, \qquad j\to \infty.
\]
Then there exists a probability measure $\sigma$ on $\T$, which is supported on a set of Lebesgue measure zero and satisfies
\[
\sum_n \abs{\widehat{\sigma}(n)}^2 w_{|n|} < \infty.
\]
\end{thm}
Below, we shall retain the notation from Subsection 3.1.
\begin{proof}[Proof of \thref{THM:singl2w}] The assumption asserts that there exists a number $0<\varepsilon <1$ and positive integers $(N_j)_j$, such that 
\[
W_j := \sup_{n: |n-N_j|\leq \varepsilon N_j} w_n \to 0, \qquad j\to \infty.
\]
By means of passing to a subsequence and of $(N_j)_j$, we may assume that $\kappa :=\inf_{j} N_{j+1}/N_j \geq 2/\varepsilon$. Passing to yet another subsequence, if necessary, we may also assume that $(W_j)_j$ is decreasing and satisfies
\[
\sum_j W_{j} 2^j < \infty.
\]
Pick an arbitrary sequence of positive numbers $(a_j)_j$ which satisfies the properties
\[
0<a_j <1, \qquad a_j \to 0, \qquad \sum_j a^2_j = +\infty.
\]
Now consider the Riesz product associated to these parameters $(N_j)_j$ and $(a_j)_j$, defined as
\[
d\sigma(\zeta) = \lim_{N\to \infty} \prod_{j=1}^N \left(1+a_j \Re(\zeta^{N_j}) \right)dm(\zeta), \qquad \zeta \in \T,
\]
where the limit is interpreted in the sense of weak-star convergence of measures. The assumptions on $(a_j)_j$ and Theorem 7.7 from Ch. V in \cite{zygmundtrigseries} ensure that $\sigma$ is probability measure supported on a set of Lebesgue measure zero. Now recall that since $\kappa > 2/\varepsilon$, the property in \eqref{EQ:Lambdaj} ensures that 
\[
\Lambda_j \subseteq \{ n: |n-N_j|\leq \varepsilon N_j \}, \qquad j=1,2, 3 \dots.
\]
Using this and the plain fact that $0<a_j<1$, the following crude estimate holds:
\[
\sum_{n\in \Lambda_{j+1}} \abs{\widehat{\sigma}(n)}^2 w_n \leq W_{j+1} \sum_{n\in \Lambda_{j+1}} \abs{\widehat{\sigma}(n)}^2 \leq W_{j+1} a^2_{j+1} \prod_{k=1}^j (1+ a^2_k) \leq W_{j+1} 2^{j+1}, \qquad j=0,1,2,\dots.
\]
Summing in $j\geq 1$, we arrive at
\[
\sum_{n > N_1} \abs{\widehat{\sigma}(n)}^2 w_{|n|} \lesssim \sum_{j\geq 1} W_{j} 2^j <\infty.
\]
\end{proof}

\bibliographystyle{siam}
\bibliography{mybib}

\Addresses

\end{document}